\documentclass{amsart}
\usepackage{amsthm}
\usepackage{amssymb}
\usepackage{upref}
\usepackage{amscd}

\numberwithin{equation}{section}
\newtheorem{main}{Theorem}

\newtheorem{theorem}[equation]{Theorem}
\newtheorem{lemma}[equation]{Lemma}
\newtheorem{proposition}[equation]{Proposition}
\newtheorem{corollary}[equation]{Corollary}

\theoremstyle{definition}
\newtheorem{definition}[equation]{Definition}

\DeclareMathOperator{\Spec}{Spec}

\DeclareMathOperator{\id}{id}
\DeclareMathOperator{\dom}{dom}
\DeclareMathOperator{\codom}{codom}
\DeclareMathOperator{\inv}{inv}

\newcommand{\cat}[1]{\mathcal{#1}}
\newcommand{\Ring}{\mathbf{Rings}}
\newcommand{\Group}{\mathbf{Groups}}
\newcommand{\Gpd}{\mathbf{Gpds}}
\newcommand{\Aff}{\mathbf{Aff}}

\newcommand{\Set}{\mathbf{Set}}
\newcommand{\Sh}{\mathbf{Sh}}
\newcommand{\site}{\Aff _{\topology }}
\newcommand{\topology}{\mathcal{T}}

\newcommand{\sheaves}{\Sh ^{\topology }}
\newcommand{\qcsheaves}{\Sh ^{\textup{qc}}}

\newcommand{\Z}{\mathbb{Z}}
\newcommand{\Fp}{\mathbb{F}_{p}}

\newcommand{\leftscript}[2]{\!\,_{#2}#1}
\newcommand{\bimod}[3]{\!\,_{#2} #1 _{#3}}

\newcommand{\mathcolon}{\colon\,}

\newcommand{\usc}{\textup{;}}

\begin{document}
 
\title{Morita theory for Hopf algebroids and presheaves of groupoids}

\date{\today}

\author{Mark Hovey}
\address{Department of Mathematics \\ Wesleyan University
\\ Middletown, CT 06459}
\thanks{The author was supported in part by NSF grant DMS 99-70978.}
\email{hovey@member.ams.org}


\begin{abstract}
Comodules over Hopf algebroids are of central importance in algebraic
topology.  It is well-known that a Hopf algebroid is the same thing as a
presheaf of groupoids on $\Aff $, the opposite category of commutative
rings.  We show in this paper that a comodule is the same thing as a
quasi-coherent sheaf over this presheaf of groupoids.  We prove the
general theorem that internal equivalences of presheaves of groupoids
with respect to a Grothendieck topology $\topology $ on $\Aff $ give
rise to equivalences of categories of sheaves in that topology.
We then show using faithfully flat descent that an internal equivalence
in the flat topology gives rise to an equivalence of categories of
quasi-coherent sheaves.  The corresponding statement for Hopf algebroids
is that weakly equivalent Hopf algebroids have equivalent categories of
comodules.  We apply this to formal group laws, where we get
considerable generalizations of the Miller-Ravenel~\cite{miller-ravenel}
and Hovey-Sadofsky~\cite{hovey-sadofsky-picard} change of rings
theorems in algebraic topology.  
\end{abstract}

\maketitle

\section*{Introduction}

A commutative Hopf algebra is a (commutative) ring $A$ together with a
lift of the functor $\Spec A\mathcolon \Ring \xrightarrow{}\Set $ to a
functor $\Ring \xrightarrow{}\Group $.  Here $\Ring $ is the category of
commutative rings with unity, $\Set $ is the category of sets, $\Group $
is the category of groups, and $(\Spec A)(R)=\Ring (A,R)$.  So a Hopf
algebra is the same thing as an affine algebraic group scheme, or a
representable presheaf of groups on $\Aff $, the opposite category of
$\Ring $.  In the same way, a Hopf algebroid $(A, \Gamma )$ is an affine
algebraic groupoid scheme, or a representable presheaf of groupoids
$(\Spec A,\Spec \Gamma )$ on $\Aff $.  Here, given a ring $R$, $\Spec
A(R)$ is the set of objects of the groupoid corresponding to $R$, and
$\Spec \Gamma (R)$ is the set of morphisms of that groupoid.  

Hopf algebroids are very important in algebraic topology, because for
many important homology theories $E$, the ring of stable co-operations
$E_{*}E$ is a (graded) Hopf algebroid over $E_{*}$ but not a Hopf
algebra.  In particular, this is true for complex cobordism $MU$ and
complex $K$-theory.  In this case, $E_{*}X$ is a (graded) comodule over
the Hopf algebroid $E_{*}E$.

Of course, not all schemes are affine.  One of the essential
contributions of Grothendieck was the realization that it is necessary
to study all schemes even if one is only interested in affine schemes.
In the same way, to understand Hopf algebroids, one should study more
general groupoid schemes.  

One of the difficulties is that the standard approach to schemes,
involving covers by open affine subschemes, is not the right one for the
algebraic topology setting.  Instead, it is better to use the functorial
approach hinted at above in our definition of $\Spec A$.  This approach
is well-known in algebraic geometry~\cite{demazure-gabriel}.  It was
introduced to algebraic topology by Hopkins and Neil Strickland.
Strickland has written an excellent exposition of this point of view
in~\cite{strickland-formal-schemes}.  In this approach, we study
arbitrary presheaves of sets (or groupoids) on $\Aff $.  

Demazure and Gabriel~\cite{demazure-gabriel} show that the category of
$A$-modules is equivalent to the category of quasi-coherent sheaves over
the presheaf of sets $\Spec A$ on $\Aff $.  Our first goal in this paper
is to extend this theorem as follows.  Let $\topology $ denote a
Grothendieck topology on $\Aff $, and let $\site $ denote the resulting
site (we put a cardinality restriction on rings to make $\Aff $ a small
category).  Given a presheaf of groupoids $(X_{0},X_{1})$ on $\Aff $, we
define the category $\sheaves _{(X_{0},X_{1})}$ of sheaves over
$(X_{0},X_{1})$ with respect to $\topology $ and we define the category
$\qcsheaves _{(X_{0},X_{1})}$ of quasi-coherent sheaves over
$(X_{0},X_{1})$.  Our first main result is then the following theorem,
proved as Theorem~\ref{thm-sheaves-hopf}.  

\begin{main}\label{main-A}
Suppose $(A,\Gamma )$ is a Hopf algebroid.  Then there is an equivalence
of categories between $\Gamma $-comodules and quasi-coherent sheaves
over $(\Spec A,\Spec \Gamma )$.
\end{main}

There is a natural notion of an internal equivalence of presheaves of
groupoids on $\site $, studied by Joyal and 
Tierney~\cite{joyal-tierney} and other authors as well.  A map $\Phi
\mathcolon (X_{0},X_{1})\xrightarrow{}(Y_{0},Y_{1})$ of presheaves of
groupoids is an internal equivalence with respect to $\topology $ if
$\Phi (R)$ is fully faithful for all $R$ and if $\Phi $ is essentially
surjective in a sheaf-theoretic sense, related to $\topology $.  This is
really the natural notion of internal equivalence for \emph{sheaves} of
groupoids on $\site $; there is a more general notion appropriate for
presheaves, introduced by Hollander~\cite{hollander}, but we do not need it.  

Our second main result is that the category of sheaves is invariant
under internal equivalence.  The following theorem is proved as
Theorem~\ref{thm-sheaves-internal}.  

\begin{main}\label{main-B}
Suppose $\Phi \mathcolon (X_{0},X_{1})\xrightarrow{}(Y_{0},Y_{1})$ is an
internal equivalence of presheaves of groupoids on $\site $.  Then $\Phi
^{*}\mathcolon \sheaves _{(Y_{0},Y_{1})}\xrightarrow{}\sheaves
_{(X_{0},X_{1})}$ is an equivalence of categories.  
\end{main}

What we really care about is the category of quasi-coherent sheaves.
Faithfully flat descent shows that a quasi-coherent sheaf is a sheaf in
the flat topology on $\Aff $.  This is often called the fpqc topology;
in it, a cover of a ring $R$ is a finite family $\{R\xrightarrow{}S_{i}
\}$ of flat extensions of $R$ such that $\prod S_{i}$ is faithfully flat
over $R$.  A strengthening of faithfully flat descent then leads to the
following theorem, proved as Theorem~\ref{thm-qc-ess-surj}. 

\begin{main}\label{main-C}
Suppose $\Phi \mathcolon (X_{0},X_{1})\xrightarrow{}(Y_{0},Y_{1})$ is an
internal equivalence of presheaves of groupoids on $\site $, where
$\topology $ is the flat topology.  Then $\Phi ^{*}\mathcolon \qcsheaves
_{(Y_{0},Y_{1})}\xrightarrow{}\qcsheaves _{(X_{0},X_{1})}$ is an
equivalence of categories.
\end{main}

In order to apply this theorem to Hopf algebroids, we need to
characterize those maps of Hopf algebroids that induce internal
equivalences in the flat topology of the corresponding presheaves of
groupoids.  The following theorem is proved as Theorem~\ref{thm-Hopf}.  

\begin{main}\label{main-D}
Suppose $f=(f_{0},f_{1})\mathcolon (A,\Gamma )\xrightarrow{}(B,\Sigma )$
is a map of Hopf algebroids.  Then $f^{*}\mathcolon (\Spec B,\Spec
\Sigma )\xrightarrow{}(\Spec A,\Spec \Gamma )$ is a internal equivalence
in the flat topology if and only if
\[
\eta _{L}\otimes f_{1}\otimes \eta _{R}\mathcolon B\otimes _{A}\Gamma
\otimes _{A}B \xrightarrow{} \Sigma 
\]
is an isomorphism and there is a ring map $g\mathcolon B\otimes
_{A}\Gamma \xrightarrow{}C$ such that $g(f_{0}\otimes \eta _{R})$
exhibits $C$ as a faithfully flat extension of $A$.
\end{main}

This condition has appeared before, in~\cite{hopkins-hopf}
and~\cite{hovey-sadofsky-picard}.  We point out that if we used the more
general notion of internal equivalence mentioned above,
Theorem~\ref{main-D} would remain unchanged, since $\Spec A$ is already
a sheaf in the flat topology by faithfully flat descent.  

Finally, we apply our results to the Hopf algebroids relevant to
algebraic topology.  The following theorem is proved as
Theorem~\ref{thm-awesome} (and the terminology is defined in
Section~\ref{sec-formal}).  

\begin{main}\label{main-E}
Fix a prime $p$ and an integer $n>0$.  Let $(A,\Gamma )$ denote the Hopf
algebroid $(v_{n}^{-1}BP_{*}/I_{n}, v_{n}^{-1}BP_{*}BP/I_{n})$.  Suppose
$B$ is a ring equipped with a homogeneous $p$-typical formal group law
of strict height $n$, classified by $f\mathcolon A\xrightarrow{}B$.
Then the functor that takes an $(A,\Gamma )$-comodule $M$ to $B\otimes
_{A}M$ defines an equivalence of categories from graded $(A,\Gamma
)$-comodules to graded $(B,B\otimes _{A}\Gamma \otimes
_{A}B)$-comodules.
\end{main}

As an immediate corollary, we recover a strengthening of the change of
rings theorem of~\cite{hovey-sadofsky-picard}, which itself is a
strengthening of the well-known Miller-Ravenel change of rings
theorem~\cite{miller-ravenel}.  The precise change of rings theorem is
prove is stated below.  

\begin{main}\label{main-F}
Let $p$ be a prime and $m\geq n>0$ be integers.  Suppose $M$ and $N$ are
$BP_{*}BP$-comodules such that $v_{n}$ acts isomorphically on $N$.  If
either $M$ is finitely presented, or if $N=v_{n}^{-1}N'$ where $N'$ is
finitely presented and $I_{n}$-nilpotent, then
\[
\text{Ext}^{**}_{BP_{*}BP}(M,N) \cong
\text{Ext}^{**}_{E(m)_{*}E(m)}(E(m)_{*}\otimes
_{BP_{*}}M,E(m)_{*}\otimes _{BP_{*}}N).  
\]
\end{main}

This theorem implies that the chromatic spectral sequence based on
$E(m)$ is the truncation of the chromatic spectral sequence based on
$BP$ consisting of the first $n+1$ columns, as pointed out
in~\cite[Remark~5.2]{hovey-sadofsky-picard}.  

There are several ways in which the results in this paper might be
generalized.  Most substantively, we do not recover the Morava change of
rings theorem~\cite{morava-comodules} from our result.  The Morava
change of rings theorem is about complete comodules over a complete Hopf
algebroid, so one would need to account for the topology in some way.
Secondly, our results will probably hold if we replace $\Aff $ by the
opposite category of rings in some topos, as suggested by Rick Jardine.
In fact, we already have to replace $\Aff $ by the opposite category of
graded rings in order to cope with the graded Hopf algebroids that arise
in algebraic topology.  Lastly, there is the aforementioned
generalization of the notion of internal equivalence, due to
Hollander~\cite{hollander}.  In this generalization, one would replace
``faithful'' by ``sheaf-theoretically faithful'' and ``full'' by
``sheaf-theoretically full''.  We are confident our results will hold
for this generalization, but we would not get any new examples of
equivalences of categories of comodules.  Nevertheless, this
generalization might be useful in other circumstances.

This paper arose from trying to understand comments of Mike Hopkins, and
I thank him deeply for sharing his insights.  The one-line summary of
this paper is ``The category of comodules over a Hopf algebroid only
depends on the associated stack'', and the author first heard this
summary from Hopkins.  It is certain that Hopkins has proved some of the
theorems in this paper.  As far as I know, however, Hopkins approached
these theorems by using stacks, which I have completely avoided.  In
particular, my definition of sheaves and quasi-coherent sheaves over
presheaves of groupoids is quite different from the definition I have
heard from Hopkins, though the two definitions are presumably
equivalent.

The author would also like to thank Dan Christensen and Rick Jardine,
both of whom thought that the original version of this paper, dealing as
it did with only quasi-coherent sheaves, was much too specific and must
be a corollary of a simpler, more general theorem.  

\section*{Notation}

We compile the notations and conventions we use in this paper.  All
rings are assumed commutative, and of cardinality less than some fixed
infinite cardinal $\kappa $.  $\Ring $ denotes the category of such
rings, and $\Aff $ denotes its opposite category.  We think of $\Aff $
as the category of representable functors $\Spec A \mathcolon \Ring
\xrightarrow{}\Set $, where $(\Spec A)(R)=\Ring (A,R)$.  We will also
want to consider $\Ring _{*}$, the category of graded rings (of
cardinality less than $\kappa $) that are commutative in the graded
sense, and its opposite category $\Aff _{*}$.

If $x,y\mathcolon A\xrightarrow{}R$ are ring homomorphisms, the symbol
$\bimod{R}{x}{y}$ denotes $R$ with its $A$-bimodule structure, where $A$
acts on the left through $x$ and on the right through $y$.  This is
especially useful for the tensor product; the symbol $R_{x}\otimes _{A}
\bimod{S}{y}{}$ indicates the bimodule tensor product, where $A$ acts
on the right on $R$ via $x$ and on the left on $S$ via $y$.  We use this
same notation in the graded case as well, where $x$ and $y$ are tacitly
assumed to preserve the grading and the tensor product is the graded
tensor product.  

The symbols $(A,\Gamma )$ and $(B,\Sigma )$ denote (possibly graded)
Hopf algebroids.  We follows the notation of~\cite[Appendix~1]{ravenel}
for the structure maps of a Hopf algebroid.  So we have the counit
$\epsilon \mathcolon \Gamma \xrightarrow{}A$, the left and right units
$\eta _{L},\eta _{R}\mathcolon A\xrightarrow{}\Gamma $, the diagonal
$\Delta \mathcolon \Gamma \xrightarrow{} \Gamma _{\eta _{R}}\otimes
_{A}\bimod{\Gamma }{\eta _{L}}{}$, and the conjugation $c\mathcolon
\bimod{\Gamma }{\eta _{L}}{\eta _{R}}\xrightarrow{}\bimod{\Gamma }{\eta
_{R}}{\eta _{L}}$.

Capital letters at the end of the alphabet, such as $X$, $Y$, and $Z$,
will denote functors from $\Ring $ to $\Set $, or functors from $\Ring
_{*}$ to $\Set $ in the graded case.  The symbol $Y_{f}\times
_{X} \bimod{Z}{g}{}$ will denote the pullback of the diagram $Y
\xrightarrow{f}X \xleftarrow{g} Z$.  

The symbols $(X_{0},X_{1})$ and $(Y_{0},Y_{1})$ will denote functors
from $\Ring $ (or $\Ring _{*}$) to $\Gpd $, the category of small
groupoids.  Here $X_{0}(R)$ is the object set of the groupoid
corresponding to $R$, and $X_{1}(R)$ is the morphism set of that
groupoid.  There are structure maps
\begin{align*}
\id & \mathcolon X_{0} \xrightarrow{} X_{1} \\
\dom ,\codom & \mathcolon X_{1} \xrightarrow{} X_{0} \\
\circ & \mathcolon (X_{1})_{\dom }\times _{X_{0}}
\leftscript{(X_{1})}{\codom }X_{1} \xrightarrow{} X_{1} \\
\inv & \mathcolon X_{1} \xrightarrow{} X_{1}
\end{align*}
satisfying the relations necessary to make $(X_{0}(R),X_{1}(R))$ a
groupoid.  

\section{Sheaves over functors}\label{sec-sheaves}

The object of this section is to define the notion of a sheaf of modules
$M$ over a sheaf of sets $X$ on $\Aff $.  We will generalize this
in the next section to sheaves of modules over sheaves of groupoids
$(X_{0},X_{1})$ on $\Aff $.  

We will assume given a Grothendieck topology $\topology $ on $\Aff $,
and denote the resulting site consisting of $\Aff $ together with
$\topology $ by $\site $.  For us, the two most important Grothendieck
topologies on $\Aff $ will be the trivial topology, where the only
covers are isomorphisms, and the the fpqc, or flat, topology, which will
be discussed later.

Now suppose $X\mathcolon \Ring \xrightarrow{}\Set $ is a functor.  We
think of $X$ as a presheaf of sets on $\site $.  We need to define the
category of sheaves over $X$.  We first define the overcategory $\site
/X$.  An object of $\site /X$ is a map of presheaves $\Spec
R\xrightarrow{x}X$, and the morphisms are the commutative triangles.  We
call the opposite category of $\site /X$ the category of \emph{points}
of $X$ following~\cite{strickland-formal-schemes}; it is called the
category of $X$-models in~\cite{demazure-gabriel}.  A point of $X$ is a
pair $(R,x)$, where $R$ is a ring and $x\in X(R)$, and a morphism from
$(R,x)$ to $(S,y)$ is a ring homomorphism $f\mathcolon R\xrightarrow{}S$
such that $X(f)(x)=y$.  We often abuse notation and write $f(x)$ for
$X(f)(x)$.  As an overcategory, $\site /X$ inherits the Grothendieck
topology $\topology $.  A cover of $(R,x)$ is a family
$\{(R,x)\xrightarrow{}(S_{i},x_{i}) \}$ such that
$\{R\xrightarrow{}S_{i} \}$ is a cover of $R$.  The category $\site /X$
also comes equipped with a structure presheaf $\cat{O}\mathcolon (\site
/X)^{\textup{op}}\xrightarrow{} \Ring $, where $\cat{O}(R,x)=R$.

\begin{definition}\label{defn-sheaves}
Suppose $X\mathcolon \Ring \xrightarrow{}\Set $ is a presheaf of sets on
$\site $.  Then a \emph{sheaf of modules over $X$}, often called just a
\emph{sheaf over $X$}, is a sheaf of $\cat{O}$-modules on $\site /X$.  
\end{definition}

More concretely, a sheaf $M$ is a functorial assignment of an $R$-module
$M_{x}$ to each point $(R,x)$, satisfying the sheaf condition.
Functoriality means that a map $(R,x)\xrightarrow{f}(S,y)$ induces a map
of $R$-modules $M_{x}\xrightarrow{\theta _{M}(f,x)}M_{y}$, where $M_{y}$
is thought of as an $R$-module by restriction.  We often abbreviate
$\theta (f,x)$ to $\theta (f)$.  We must have $\theta 
(gf)=\theta  (g)\circ \theta (f)$ and $\theta (1)=1$.  The sheaf
condition means that if $\{(R,x)\xrightarrow{}(S_{i},x_{i}) \}$ is a
cover, then the diagram 
\[
M_{x} \xrightarrow{} \prod_{i} M_{x_{i}} \rightrightarrows \prod _{jk}
M_{x_{jk}} 
\]
is an equalizer of $R$-modules, where $x_{jk}$ is the image of $x$ in
$X(S_{j}\otimes _{R}S_{k})$.  The maps in this diagram are all maps of
$R$-modules.  

We have an evident definition of a map of sheaves over $X$.  To be
concrete, a map $\alpha \mathcolon M\xrightarrow{}N$ of sheaves over $X$
assigns to each point $(R,x)$ of $X$ a map $\alpha _{x}\mathcolon
M_{x}\xrightarrow{}N_{x}$ of $R$-modules, natural in $(R,x)$.  This
gives us a category $\sheaves _{X}$ of sheaves over $X$.  A map of
sheaves $X\xrightarrow{\Phi }Y$ induces a functor $\Phi ^{*}\mathcolon
\sheaves _{Y} \xrightarrow{}\sheaves _{X}$. Here, if $M$ is a sheaf over
$Y$ and $(R,x)$ is a point of $X$, we define $(\Phi ^{*}M)_{x}=M_{\Phi
x}$.    

Note that all of these definitions work perfectly well in the graded
case as well.  We would have a Grothendieck topology $\topology $ on
$\Aff _{*}$, and a functor $X\mathcolon \Aff _{*}\xrightarrow{}\Set $.
A point of $X$ would be a graded ring $R$ and a point $x\in X(R)$.  A
sheaf $M$ over $X$ would be as assignment of a graded $R$-module $M_{x}$
to each point $(R,x)$ of $X(R)$, satisfying the functoriality and sheaf
conditions.  

We now consider quasi-coherent sheaves.  We only need quasi-coherent
sheaves in the trivial topology, so we will stick to that case.  A
quasi-coherent sheaf is supposed to be a sheaf that acts like a free
sheaf in an appropriate sense.  The salient property of the free sheaf
$\cat{O}$ is that, if $(R,x)\xrightarrow{}(S,y)$ is a map of points,
then $\cat{O}_{y}=S\otimes _{R}\cat{O}_{x}$.  We therefore make the
following definition.

\begin{definition}\label{defn-qc}
Suppose $X\mathcolon \Ring \xrightarrow{}\Set $ is a functor.  A
\emph{quasi-coherent sheaf $M$ over $X$} is a sheaf over $X$ in the
trivial topology such that, given a map $(R,x)\xrightarrow{}(S,y)$ of
points of $X$, the adjoint $\rho _{M}(f)\mathcolon S\otimes
_{R}M_{x}\xrightarrow{}M_{y}$ of $\theta _{M}(f)$ is an isomorphism.
\end{definition}

This is the same definition given in~\cite{demazure-gabriel}
and~\cite{strickland-formal-schemes}.  We get a category
$\qcsheaves _{X}$, which is the full subcategory of sheaves over $X$ in
the trivial topology consisting of the quasi-coherent sheaves.  Given a
map $\Phi \mathcolon X\xrightarrow{}Y$ of functors, $\Phi ^{*}\mathcolon
\sheaves _{Y}\xrightarrow{}\sheaves _{X}$ restricts to define $\Phi
^{*}\mathcolon \qcsheaves _{Y}\xrightarrow{}\qcsheaves _{X}$.  

The value of this definition of quasi-coherence is shown by the
following lemma. 

\begin{lemma}\label{lem-qc-module}
Suppose $A\in \Ring $, and let $\Spec A\mathcolon \Ring
\xrightarrow{}\Set $ be the representable functor $(\Spec A)(R)=\Ring
(A,R)$.  Then the category of $A$-modules is equivalent to the category
of quasi-coherent sheaves over $\Spec A$.  The equivalence takes an
$A$-module $M$ to the quasi-coherent sheaf $\widetilde{M}$ over $\Spec
A$ defined by $\widetilde{M}_{x}=R_{x}\otimes _{A}M$ for $x\mathcolon
A\xrightarrow{}R$, and its inverse takes a quasi-coherent sheaf $N$ to
its value at $1\mathcolon A\xrightarrow{}A$.  
\end{lemma}

This lemma is due to Demazure and
Gabriel~\cite[p.~61]{demazure-gabriel-intro}, who actually show that the
category of quasi-coherent sheaves over a scheme when defined this way
agrees (up to equivalence) with the usual notion of quasi-coherent
sheaves on a scheme.  A direct proof can be found
in~\cite{strickland-formal-schemes}.

Once again, we note that Lemma~\ref{lem-qc-module} will work in the
graded case as well.  The definition of a quasi-coherent sheaf over a
functor $X\mathcolon \Ring _{*}\xrightarrow{}\Set $ is similar to the
ungraded case, and the same argument used to prove
Lemma~\ref{lem-qc-module} shows that, if $A$ is a graded ring, the
category of quasi-coherent sheaves over $\Spec A$ (now defined by
$(\Spec A)(R)=\Ring _{*}(A,R)$) is equivalent to the category of graded
$A$-modules.

It will be useful later to note that, if $f\mathcolon A\xrightarrow{}B$
is a ring homomorphism and $\Spec f\mathcolon \Spec B\xrightarrow{}\Spec A$
is the corresponding map of functors, then the induced map $(\Spec
f)^{*}\mathcolon \qcsheaves _{\Spec A}\xrightarrow{}\qcsheaves _{\Spec
B}$ takes the $A$-module $M$ to the $B$-module $B\otimes _{A}M$.  

\section{Sheaves over groupoid functors}\label{sheaves-gr}

The object of this section is to prove Theorem~\ref{main-A}, showing
that a comodule over a Hopf algebroid is a special case of the more
general notion of a quasi-coherent sheaf over a presheaf of groupoids.
This will require us to define the notion of a sheaf $M$ of modules over
a presheaf of groupoids $(X_{0},X_{1})$ on $\site $.

We will consider a presheaf of groupoids $(X_{0},X_{1})$ on $\site $.
This means that $X_{0}$ and $X_{1}$ are presheaves of sets on $\site $,
and that $(X_{0}(R),X_{1}(R))$ is a groupoid for all $R$, naturally in
$R$.  So we have structure maps as defined in the notation section.  A
presheaf of groupoids $(X_{0},X_{1})$ is called a \emph{sheaf of
groupoids} when $X_{0}$ and $X_{1}$ are sheaves of sets on $\site $; we
would be happy to assume our presheaves of groupoids are in fact sheaves
of groupoids, but that assumption is unnecessary.  

Sheaves of groupoids have been much studied in the literature; a stack
is a special kind of sheaf of groupoids, and stacks are essential in
modern algebraic geometry~\cite{faltings-chai}.  The homotopy theory of
sheaves of groupoids has been studied by Joyal and
Tierney~\cite{joyal-tierney}, Jardine~\cite{jardine-stacks}, and
Hollander~\cite{hollander}.  

\begin{definition}\label{defn-sheaf-groupoid}
Suppose $(X_{0},X_{1})$ is a presheaf of groupoids on $\site $.  A
\emph{sheaf over $(X_{0},X_{1})$} is a sheaf $M$ over $X_{0}$ together
with an isomorphism $\psi \mathcolon \dom ^{*}M\xrightarrow{}\codom
^{*}M$ of sheaves over $X_{1}$ satisfying the cocycle condition.  To
explain the cocycle condition, note that, if $\alpha $ is a morphism of
$X_{1}(R)$, $\psi _{\alpha }$ is an isomorphism of $R$-modules $\psi
_{\alpha }\mathcolon M_{\dom \alpha }\xrightarrow{}M_{\codom \alpha }$.
The cocycle condition says that if $\beta $ and $\alpha $ are composable
morphisms, then $\psi _{\beta \alpha }=\psi _{\beta }\circ \psi _{\alpha
}$.  A \emph{quasi-coherent sheaf over $(X_{0},X_{1})$} is a sheaf $M$
over $(X_{0},X_{1})$ in the trivial topology such that $M$ is
quasi-coherent as a sheaf over $X_{0}$.  
\end{definition}

We also get a notion of a map $\tau \mathcolon M\xrightarrow{}N$ of
sheaves over $(X_{0},X_{1})$.  Such a map is a map of sheaves over
$X_{0}$ such that the diagram 
\[
\begin{CD}
M_{\dom \alpha } @>\psi ^{M}_{\alpha } >> M_{\codom \alpha } \\
@V\tau _{\dom \alpha } VV @VV\tau _{\codom \alpha }V \\
N_{\dom \alpha } @>>\psi ^{N}_{\alpha }> N_{\codom \alpha }
\end{CD}
\]
commutes for all points $(R,\alpha )$ of $X_{1}(R)$.  We then get 
categories $\sheaves _{(X_{0},X_{1})}$ and $\qcsheaves
_{(X_{0},X_{1})}$.  

Note that a map $\Phi \mathcolon
(X_{0},X_{1})\xrightarrow{}(Y_{0},Y_{1})$ induces a functor $\Phi
^{*}\mathcolon \sheaves _{(Y_{0},Y_{1})}\xrightarrow{} \sheaves
_{(X_{0},X_{1})}$ and $\Phi ^{*}\mathcolon \qcsheaves
_{(Y_{0},Y_{1})}\xrightarrow{}\qcsheaves _{(X_{0},X_{1})}$.  Indeed, we
define $\psi ^{\Phi ^{*}M}_{\alpha }=\psi ^{M}_{\Phi \alpha }$.  

Also note that all of the comments above work perfectly well for
presheaves of groupoids on $\Aff _{*}$.  In this case, $\psi _{\alpha
}\mathcolon M_{\dom \alpha }\xrightarrow{}M_{\codom \alpha }$ will be an
isomorphism of graded $R$-modules.  

As originally noted by Haynes Miller, a Hopf
algebroid~\cite[Appendix~1]{ravenel} is just a pair of commutative rings
$(A, \Gamma )$ such that $(\Spec A,\Spec \Gamma )$ is a sheaf of
groupoids (in the trivial topology).  The structure maps of a Hopf
algebroid (listed in the notation section) are therefore dual to the
structure maps of a presheaf of groupoids; for example, the diagonal
$\Delta \mathcolon \Gamma \xrightarrow{}\Gamma _{\eta _{R}}\otimes
_{A}\bimod{\Gamma }{\eta _{L}}{}$ is dual to the composition map
$(X_{1})_{\dom }\times _{X_{0}}\bimod{X}{\codom }{1}$.

It is useful to recall the composition in the groupoid $(\Spec A,\Spec
\Gamma )(R)$ from this point of view.  Suppose $\beta ,\alpha \mathcolon
\Gamma \xrightarrow{}R$ are ring homomorphisms with $\alpha \eta
_{L}=x$, $\alpha \eta _{R} = \beta \eta _{L}=y$, and $\beta \eta
_{R}=z$, so that $\alpha $ is a morphism from $x$ to $y$ and $\beta $ is
a morphism from $y$ to $z$.  The composition $\beta \circ \alpha
\mathcolon \Gamma \xrightarrow{}R$ is defined to be the composite
\[
\Gamma \xrightarrow{\Delta } \bimod{\Gamma }{\eta _{L} }{\eta _{R}}
\otimes _{A} \bimod{\Gamma }{\eta _{L} }{\eta _{R}} \xrightarrow{\alpha
\otimes \beta } \bimod{R}{x}{y} \otimes _{A} \bimod{R}{y}{z}
\xrightarrow{\mu } \bimod{R}{x}{z}.
\]

Just as a quasi-coherent sheaf over $\Spec A$ is the same thing as a
module over $A$, so a quasi-coherent sheaf over $(\Spec A,\Spec \Gamma
)$ is the same thing as a comodule over $(A,\Gamma )$.  The following
theorem is Theorem~\ref{main-A} of the introduction.  

\begin{theorem}\label{thm-sheaves-hopf}
Suppose $(A,\Gamma )$ is a Hopf algebroid.  Then there is an equivalence
of categories between $\Gamma $-comodules and quasi-coherent sheaves
over $(\Spec A,\Spec \Gamma )$.
\end{theorem}

This theorem will also hold in the graded context: if $(A,\Gamma )$ is a
graded Hopf algebroid, then the category of graded $\Gamma $-comodules
is equivalent to the category of quasi-coherent sheaves over the
presheaf of groupoids $(\Spec A,\Spec \Gamma )$ on $\Aff _{*}$.  The
proof is the same as the proof below.    

\begin{proof}
We first construct a functor from quasi-coherent sheaves over
$(\Spec A,\Spec \Gamma )$ to $(A,\Gamma )$-comodules.  Suppose that
$\widetilde{M}$ is a quasi-coherent sheaf over $(\Spec A,\Spec \Gamma
)$.  Then $\widetilde{M}$ is in particular a quasi-coherent sheaf over
$\Spec A$, so corresponds to an $A$-module $M$.  Then if $\alpha
\mathcolon \Gamma \xrightarrow{}R$ is a point of $\Spec \Gamma $ defined
over $R$, with $\alpha \eta _{L}=x$ and $\alpha \eta _{R}=y$,
\[
(\dom ^{*}\widetilde{M})_{\alpha }=R_{x}\otimes _{A} M \text{ and }
(\codom ^{*}\widetilde{M})_{\alpha }=R_{y}\otimes _{A}M.
\]
Let us denote by $\widetilde{\psi }$ the isomorphism of sheaves $\dom
^{*}\widetilde{M}\xrightarrow{}\codom ^{*}\widetilde{M}$.  Then,
$\widetilde{\psi }$ defines an isomorphism
\[
\widetilde{\psi }_{\alpha }\mathcolon R_{x}\otimes _{A} M \xrightarrow{}
R_{y} \otimes _{A} M.
\]
of $R$-modules.  Taking $\alpha $ to be the identity map $1$ of $\Gamma
$, we define $\psi \mathcolon M\xrightarrow{}\Gamma _{\eta _{R}}\otimes
_{A}M$ to be the composite
\[
M = A\otimes _{A}M \xrightarrow{\eta _{L}\otimes 1} \Gamma_{\eta _{L}}
\otimes _{A} M \xrightarrow{\widetilde{\psi }_{1}} \Gamma _{\eta _{R}}
\otimes _{A} M.
\]

We must show that $\psi $ is counital and coassociative.  Note first
that $\epsilon \mathcolon \Gamma \xrightarrow{}A$, thought of as a
morphism in the groupoid $(\Spec A,\Spec \Gamma )(A)$, is the identity
morphism of the object $1_{A}\mathcolon A\xrightarrow{}A$, and so in
particular is idempotent.  The cocycle condition implies that
$\widetilde{\psi }_{\epsilon }$ is also idempotent, and since it is an
isomorphism, it follows that $\widetilde{\psi }_{\epsilon }$ is the
identity of $M$.  Now, $\epsilon $ defines a map from the point $(\Gamma
,1)$ to the point $(A, \epsilon )$ of $\Spec \Gamma $.  Since
$\widetilde{\psi }$ is a map of sheaves over $\Spec \Gamma $, we
conclude that
\[
1\otimes \widetilde{\psi }_{1} \mathcolon A\otimes _{\Gamma }(\Gamma
_{\eta _{L} }\otimes _{A} M) \xrightarrow{} A\otimes _{\Gamma }(\Gamma
_{\eta _{R}} \otimes _{A} M)
\]
is the identity map.  From this it follows easily that $\psi $ is
counital.  

To see that $\psi $ is coassociative, let $\alpha \mathcolon \Gamma
\xrightarrow{}\Gamma \otimes _{A}\Gamma $ denote the map that takes
$t$ to $t \otimes 1$.  Let $\beta $ denote the map that takes
$t$ to $1\otimes t$.  Then we have 
\[
y\eta _{R}(a) = \eta _{R}a\otimes 1 = 1 \otimes \eta _{L}a =x\eta
_{L}(a),
\]
and so $\beta \circ \alpha $ makes sense.  A calculation shows that
$\beta \circ \alpha =\Delta $, the diagonal map.  If $(R,\gamma )$ is an
arbitrary point of $\Spec \Gamma $ with $\gamma \eta _{L}=x$ and $\gamma
\eta _{R}=y$, there is a map from $(\Gamma ,1)$ to $(R,\gamma )$.  Since
$\widetilde{\psi }$ is a map of sheaves, we find that $\widetilde{\psi
}_{\gamma }$ is the composite
\[
R_{x}\otimes _{A} M \cong R_{\gamma }\otimes _{\Gamma }\Gamma _{\eta
_{L} }\otimes _{A} M \xrightarrow{1\otimes \widetilde{\psi }_{1}}
R_{\gamma } \otimes _{\Gamma } \Gamma _{\eta _{R}}\otimes _{A} M \cong
R_{y}\otimes _{A} M.
\]
This description allows us to compute $\widetilde{\psi }_{\beta }$ and
$\widetilde{\psi }_{\alpha }$, and so also their composite.  We find
that $\widetilde{\psi }_{\beta }\circ \widetilde{\psi }_{\alpha }$ takes
$1\otimes 1\otimes m$ to $(1\otimes \psi )\psi (m)$.  Similarly
$\widetilde{\psi }_{\Delta }$ takes $1\otimes 1\otimes m$ to $(\Delta
\otimes 1)\psi (m)$.  The cocycle condition forces these to be equal,
and so $\psi $ is coassociative.

We have now constructed a comodule $M$ associated to any quasi-coherent
sheaf $\widetilde{M}$ over $(\Spec A,\Spec \Gamma )$.  We leave to the
reader the striaghtforward check that this is functorial.  

Our next goal is to construct a functor from $(A,\Gamma )$-comodules to
quasi-coherent sheaves over $(\Spec A,\Spec \Gamma )$.  Suppose $M$ is a
$\Gamma $-comodule with structure map $\psi \mathcolon
M\xrightarrow{}\Gamma_{\eta _{R}} \otimes _{A}M$.  Then, in particular,
$M$ is an $A$-module, so there is an associated quasi-coherent sheaf
$\widetilde{M}$ over $\Spec A$, defined by $\widetilde{M}_{x}=R_{x}\otimes
_{A}M$, where $x\mathcolon A\xrightarrow{}R$ is a ring homomorphism.
Given a point $\alpha \mathcolon \Gamma \xrightarrow{}R$ of $\Spec
\Gamma $ with $\alpha \eta _{L}=x$ and $\alpha \eta _{R}=y$, we have 
\[
(\dom ^{*}\widetilde{M})_{x}=R_{x}\otimes _{A}M \text{ and } (\codom
^{*}\widetilde{M})_{x}=R_{y}\otimes _{A}M
\]
We define $\widetilde{\psi }\mathcolon \dom ^{*}\widetilde{M}
\xrightarrow{}\codom ^{*}\widetilde{M}$ by letting $\widetilde{\psi
}_{\alpha }$ be the composite
\[
R_{x}\otimes _{A}M \xrightarrow{1\otimes \psi } R_{x} \otimes _{A}
\bimod{\Gamma }{\eta _{L} }{\eta _{R}} \otimes _{A} M
\xrightarrow{1\otimes \alpha \otimes 1} R_{x} \otimes _{A}
\bimod{R}{x}{y} \otimes _{A} M \xrightarrow{\mu \otimes 1} R_{y}\otimes
_{A}M.
\]
We leave to the reader the check that $\widetilde{\psi }$ is a map of
sheaves.  

It remains to show that $\widetilde{\psi }$ satisfies the cocycle
condition and is an isomorphism.  We begin with the cocycle condition.
Suppose that $\alpha , \beta \mathcolon \Gamma \xrightarrow{}R$ are ring
homomorphisms with $\alpha \eta _{L}=x$, $\alpha \eta _{R}=\beta \eta
_{L}=y$, and $\beta \eta _{R}=z$.  Consider the following commutative
diagram, in which all tensor products that occur are taken over $A$, and
$\Gamma =\bimod{\Gamma }{\eta _{L} }{\eta _{R}}$.  

\scriptsize
\[
\begin{CD}
R_{x} \otimes M @>1\otimes \psi >> R_{x} \otimes \Gamma \otimes M
@>1\otimes \alpha \otimes 1 >> R_{x} \otimes \bimod{R}{x}{y} \otimes M
@>\mu \otimes 1 >> R_{y} \otimes  M \\
@V1\otimes \psi VV @V1\otimes 1\otimes \psi VV @V1\otimes 1\otimes \psi
VV @V 1\otimes 1\otimes \psi VV \\
R_{x} \otimes \Gamma \otimes M @>1\otimes \Delta \otimes 1 >> R_{x}
\otimes \Gamma \otimes \Gamma \otimes M @> 1\otimes \alpha \otimes
1\otimes 1>> R_{x} \otimes \bimod{R}{x}{y} \otimes \Gamma \otimes M
@>\mu
\otimes 1\otimes 1 >> R_{y} \otimes  \Gamma \otimes  M \\
@. @. @V1\otimes 1\otimes \beta \otimes 1 VV @V1\otimes \beta  \otimes 1 VV \\
\ @. \ @. R_{x} \otimes \bimod{R}{x}{y} \otimes
\bimod{R}{y}{z} \otimes M @>\mu \otimes 1\otimes 1>> R_{y} \otimes
\bimod{R}{y}{z} \otimes M \\ 
@. @. @V1\otimes \mu \otimes 1 VV @V\mu \otimes 1VV \\
\ @. \ @. R_{x} \otimes \bimod{R}{x}{z} \otimes M @>\mu
\otimes 1>> R_{z} \otimes M
\end{CD}
\]
\normalsize 

The outer clockwise composite in this diagram is
$\widetilde{\psi }_{\beta }\circ \widetilde{\psi }_{\alpha }$, and the
outer counterclockwise composite is $\widetilde{\psi }_{\beta \circ
\alpha }$, using the description of $\beta \circ \alpha $ given above.
Thus $\widetilde{\psi }$ satisfies the cocycle condition.

We must still show that $\widetilde{\psi }_{\alpha }$ is an isomorphism
for all $\alpha \mathcolon \Gamma \xrightarrow{}R$.  Since
$\widetilde{\psi }$ satisfies the cocycle condition and $\alpha $ is
itself an isomorphism, it suffices to show that $\widetilde{\psi
}_{1_{x}}$ is an isomorphism, where $1_{x}$ is the identity morphism of
$x\mathcolon A\xrightarrow{}R$.  That is, $1_{x}$ is the composite 
\[
\Gamma \xrightarrow{\epsilon } A \xrightarrow{x} R.
\]
But one can check, using the fact that $\psi $ is counital, that
$\widetilde{\psi }_{1_{x}}$ is the identity of $R_{x}\otimes _{A}M$.  
This completes the proof that $\widetilde{M}$ is a quasi-coherent sheaf
over $(\Spec A,\Spec \Gamma )$.  We leave to the reader the check that
it is functorial in $M$.  

We also leave to the reader the check that these constructions define
inverse equivalences of categories.
\end{proof}

Maps of Hopf algebroids $(f_{0},f_{1})\mathcolon (A, \Gamma
)\xrightarrow{}(B,\Sigma )$ are defined
in~\cite[Definition~A1.1.7]{ravenel}; they are, of course, maps such
that $\Phi =(\Spec f_{0},\Spec f_{1})$ is a map of sheaves of
groupoids.  According to Theorem~\ref{thm-sheaves-hopf}, $(f_{0},f_{1})$
will induce a map $\Phi ^{*}$ from $(A, \Gamma )$-comodules to
$(B,\Sigma )$-comodules.  This maps takes the $\Gamma $-comodule $M$ to
$B\otimes _{A}M$.  In order to define the structure map of
$B\otimes _{A} M$, recall from~\cite[Definition~A1.1.7]{ravenel} that the
definition of a map of Hopf algebroids requires 
\[
\eta _{L}f_{0}=x=f_{1}\eta _{L} \text{ and } \eta _{R}f_{0}=y=f_{1}\eta
_{R}. 
\]
We then define the structure map of $B\otimes _{A}M$ to be the composite
\begin{align*}
B_{f_{0}}\otimes _{A} M & \xrightarrow{1\otimes \psi } B \otimes _{A}
\bimod{\Gamma }{\eta _{L}}{\eta _{R}} 
\otimes _{A} M  \xrightarrow{\eta _{L}\otimes f_{1}\otimes 1} \Delta _{x} 
\otimes _{A} \bimod{\Delta }{x}{y} \otimes _{A} M \\
& \xrightarrow{\mu \otimes 1} \Delta _{y}\otimes _{A} M \cong
\Delta _{\eta _{R}} \otimes _{B} (B_{f_{0}}\otimes _{A} M).  
\end{align*}

\section{Internal equivalences yield equivalences}\label{sec-equiv}

The object of this section is to prove Theorem~\ref{main-B}, showing
that if $\Phi \mathcolon (X_{0},X_{1})\xrightarrow{}(Y_{0},Y_{1})$ is an
\emph{internal equivalence} of presheaves of groupoids on $\site $, then
\[
\Phi ^{*}\mathcolon \sheaves _{(Y_{0},Y_{1})}\xrightarrow{}\sheaves
_{(X_{0},X_{1})}
\]
is an equivalence of categories.  This statement essentially says that
the category of sheaves is a homotopy-invariant construction.
 
We begin by defining an internal equivalence.  Internal equivalences are
the weak equivalences in the model structure on sheaves of groupoids
considered by Joyal and Tierney in~\cite{joyal-tierney}.

\begin{definition}\label{defn-internal}
Suppose $\Phi \mathcolon (X_{0},X_{1})\xrightarrow{}(Y_{0},Y_{1})$ is a
map of presheaves of groupoids on $\site $.  The \emph{essential image
of $\Phi $} is the subfunctor of $Y_{0}$ consisting of all points
$(R,y)$ of $Y_{0}$ such that there exists a point $(R,x)$ of $X_{0}$ and
a morphism $\alpha \in Y_{1}(R)$ from $\Phi x$ to $y$.  The
\emph{sheaf-theoretic essential image of $\Phi $} is the subfunctor of
$Y_{0}$ consisting of all points $(R,y)$ such that there exists a cover
$\{R\xrightarrow{f_{i}}S_{i} \}$ of $R$ in the topology $\topology $
such that $y_{i}=f_{i}y$ is in the essential image of $\Phi $ for all
$i$.  The map $\Phi $ is called an \emph{internal equivalence} if $\Phi
(R)$ is full and faithful for all $R$, and if the sheaf-theoretic
essential image of $\Phi $ is $Y_{0}$ itself.  
\end{definition}

For example, $\Phi $ is an internal equivalence in the trivial topology
if and only if $\Phi (R)$ is full, faithful, and essentially surjective
for all $R$, so that $\Phi (R)$ is an equivalence of groupoids for all
$R$.  

Our goal is then to prove the following theorem, which is
Theorem~\ref{main-B} of the introduction.  

\begin{theorem}\label{thm-sheaves-internal}
Suppose $\Phi \mathcolon (X_{0},X_{1})\xrightarrow{}(Y_{0},Y_{1})$ is an
internal equivalence of presheaves of groupoids on $\site $.  Then $\Phi
^{*}\mathcolon \sheaves _{(Y_{0},Y_{1})}\xrightarrow{}\sheaves
_{(X_{0},X_{1})}$ is an equivalence of categories.  
\end{theorem}

As usual, our proof of this theorem will work in the graded case as
well.  

We point out that there should be a model structure on presheaves of
groupoids extending the Joyal-Tierney model structure.  The weak
equivalences in this model structure would be the maps $\Phi $ which are
sheaf-theoretically fully faithful and whose sheaf-theoretic essential
image is all of $Y_{0}$.  Theorem~\ref{thm-sheaves-internal} should then
be a special case of the more general theorem that a weak equivalence
of presheaves of groupoids induces an equivalence of their categories of
sheaves.  We have not considered this more general case, because
$\Spec A$ is already a sheaf in the flat topology, and $\Spec A$ is our
main object of interest.  

We will prove this theorem by showing that $\Phi ^{*}$ is full,
faithful, and essentially surjective.  The proof of each such step will
be long, but divided into discrete steps very much like a diagram chase.
In general, we are trying in each case to construct something for every
point $(R,y)$ of $Y_{0}$.  So first we do it for points $(R,y)$ in the
essential image of $\Phi $.  This generally involves choosing a point
$(R,x)$ of $X_{0}$ and a morphism $\alpha \mathcolon \Phi
x\xrightarrow{}y$, so we generally have to prove that which choice one
makes is immaterial.  Then we show that every property we hope for
in the construction is true on the essential image of $\Phi $.  Next we
extend the definition to all points $(R,y)$ in the sheaf-theoretic
essential image of $\Phi $ by using a cover.  Once again, this depends
on the choice of cover, so we have to show the choice is immaterial.
For this, it is enough to show that refining the cover makes no
difference, since any two covers have a common refinement.  Finally, we
show that the properties we want are sheaf-theoretic in nature, so that
since they hold already on the essential image of $\Phi $, they also
hold on the sheaf-theoretic essential image of $\Phi $.  

\begin{proposition}\label{prop-faithful}
Suppose $\Phi \mathcolon (X_{0},X_{1})\xrightarrow{}(Y_{0},Y_{1})$ is an
map of presheaves of groupoids on $\site $ whose sheaf-theoretic essential
image is all of $Y_{0}$.  Then 
\[
\Phi ^{*}\mathcolon \sheaves _{(Y_{0},Y_{1})}\xrightarrow{}\sheaves
_{(X_{0},X_{1})}
\]
is faithful.
\end{proposition}

\begin{proof}
Suppose $\tau \mathcolon M\xrightarrow{}N$ is a map of 
sheaves on $(Y_{0},Y_{1})$ such that $\Phi ^{*}\tau =0$.  This means that
$\tau _{\Phi x}=0$ for all points $(R,x)$ of $X_{0}$.  We must show that
$\tau _{y}=0$ for all points $(R,y)$ of $Y_{0}$.  We first show that
$\tau _{y}=0$ for all $y$ in the essential image of $\Phi $.  Indeed,
suppose $\alpha $ is a morphism from $\Phi x$ to $y$.  Then, since $\tau
$ commutes with the structure map $\psi $, we get the commutative
diagram below.  
\[
\begin{CD}
M_{\Phi x} @>\psi ^{M}_{\alpha } >> M_{y} \\
@V\tau _{\Phi x} VV @VV\tau _{y}V \\
N_{\Phi x} @>>\psi ^{N}_{\alpha } > N_{y}
\end{CD}
\]
It follows that $\tau _{y}=0$.  

Now suppose $(R,y)$ is a general point of $Y_{0}$.  Since $y$ is in the
sheaf-theoretic essential image of $\Phi $,  we can choose a covering
$\{R\xrightarrow{f_{i}}S_{i} \}$ such that $y_{i}=Y_{0}(f_{i})(y)$ is
in the essential image of $\Phi $ for all $i$.  Thus $\tau _{y_{i}}=0$
for all $i$.  We then have a commutative diagram 
\[
\begin{CD}
M_{y} @>>> \prod M_{y_{i}} \\
@V\tau _{y}VV @VV\prod \tau _{y_{i}}V \\
N_{y} @>>> \prod N_{y_{i}}
\end{CD}
\]
The horizontal arrows are monomorphisms, since $M$ and $N$ are sheaves
in $\topology $, so $\tau _{y}=0$ as well.  
\end{proof}

Note that we have actually shown, more generally, that if $\tau
\mathcolon M\xrightarrow{}N$ is a morphism of sheaves over $(Y_{0},Y_{1})$
such $\Phi ^{*}\tau =0$, then $\tau $ restricted to the sheaf-theoretic
essential image of $\Phi $ is also $0$.  

\begin{proposition}\label{prop-full}
Suppose $\Phi \mathcolon (X_{0},X_{1})\xrightarrow{}(Y_{0},Y_{1})$ is an
map of presheaves of groupoids on $\site $ whose sheaf-theoretic essential
image is all of $Y_{0}$ and such that $\Phi (R)$ is full for all $R$.
Then $\Phi ^{*}\mathcolon \sheaves
_{(Y_{0},Y_{1})}\xrightarrow{}\sheaves _{(X_{0},X_{1})}$ is full.
\end{proposition}

\begin{proof}
Suppose we have a map $\tau \mathcolon \Phi ^{*}M\xrightarrow{}\Phi
^{*}N$.  This means we have maps $\tau _{x}\mathcolon M_{\Phi
x}\xrightarrow{}N_{\Phi x}$ for all points $(R,x)$ of $X_{0}$.  We need
to construct maps $\sigma _{y}\mathcolon M_{y}\xrightarrow{}N_{y}$ for
all points $(R,y)$ of $Y_{0}$ such that $\sigma _{\Phi x}=\tau _{x}$.
Suppose first that $y$ is in the essential image of $\Phi $, so that
there is a morphism $\alpha $ from $\Phi x$ to $y$ for some point
$(R,x)$ of $X_{0}$.  If $\sigma $ were to exist, then we would have the
commutative diagram below,
\[
\begin{CD}
M_{\Phi x} @>\psi ^{M}_{\alpha }>> M_{y} \\
@V\tau _{x} VV @VV\sigma _{y}V \\
N_{\Phi x} @>>\psi ^{N}_{\alpha }> N_{y}
\end{CD}
\]
so we define $\sigma _{y}=\psi ^{N}_{\alpha }\tau _{x}(\psi ^{M}_{\alpha
})^{-1}$.

We claim that this definition of $\sigma _{y}$ is independent of the
choice of $\alpha $.  Indeed, suppose $\beta \in Y_{1}(R)$ is a morphism
from $\Phi x'$ to $y$.  Then $\beta ^{-1}\alpha $ is a morphism from
$\Phi x$ to $\Phi x'$, and so, since $\Phi $ is full, there is a
morphism $\gamma \in X_{1}(R)$ from $x$ to $x'$ such that $\Phi \gamma
=\beta ^{-1}\alpha $.  Since $\tau $ is a map of sheaves, $\tau
_{x'}\psi ^{M}_{\Phi \gamma }=\psi ^{N}_{\Phi \gamma }\tau _{x}$.  On
the other hand, by the cocycle condition we have $\psi _{\Phi \gamma
}=(\psi _{\beta })^{-1}\psi _{\alpha }$.  Combining these two equations
gives
\[
\psi ^{N}_{\alpha }\tau _{x}(\psi ^{M}_{\alpha })^{-1} = \psi
^{N}_{\beta }\tau _{x'}(\psi ^{M}_{\beta })^{-1},
\]
so $\sigma _{y}$ is independent of the choice of $\alpha $.  In
particular, if $y=\Phi x$, we can take $\alpha $ to be the identity map
of $\Phi x$.  The cocycle condition implies that $\psi ^{M}_{\alpha }$
and $\psi ^{N}_{\alpha }$ are identity maps, and so $\sigma _{\Phi
x}=\tau _{x}$.

We now show that $\sigma $ commutes with the structure maps of $M$ and
$N$ on the essential image of $\Phi $.  Suppose that
$(R,y)\xrightarrow{f}(S,y')$ is a map of points of $Y_{0}$, and that $y$
is in the essential image of $\Phi $.  Choose a morphism $\alpha $ from
$\Phi x$ to $y$ for some point $(R,x)$ of $X_{0}$.  Let $\alpha
'=Y_{1}(f)(\alpha )$, so that $\alpha '$ is a morphism from $\Phi x'$ to
$y'$, where $x'=X_{0}(f)(x)$.  Since $\tau $ is a map of sheaves, we get
the commutative square below.
\[
\begin{CD}
M_{\Phi x} @>\tau _{x} >> N_{\Phi x} \\
@V\theta ^{M}(f,\Phi x)VV @VV\theta ^{N}(f,\Phi x)V \\
M_{\Phi x'} @>>\tau _{x'}> N_{\Phi x'}
\end{CD}
\]
We would like to know that the square below is commutative. 
\[
\begin{CD}
M_{y} @>\sigma _{y} >> N_{y} \\
@V\theta ^{M}(f,y)VV @VV\theta ^{N}(f,y)V \\
M_{y'} @>>\sigma _{hy}>> N_{y'}.
\end{CD}
\]
We claim that is an isomorphism from the top square to the bottom
square, and so the bottom square must be commutative.  Indeed, in the
upper left corner this isomorphism is $\psi ^{M}_{\alpha }$, in
the upper right corner it is $\psi ^{N}_{\alpha }$, in the
lower left corner it is $\psi ^{M}_{\alpha '}$, and in the lower right
corner it is $\psi ^{N}_{\alpha '}$.  All the required diagrams commute
to make this a map of squares.  This uses the fact that $\psi ^{M}$ and
$\psi ^{N}$ are maps of sheaves and the well-definedness of $\sigma $.  

We now check that $\sigma $ commutes with $\psi $, on the essential
image of $\Phi $.  Suppose we have a morphism $\beta \mathcolon
y\xrightarrow{}y'$ in $(Y_{0}(R), Y_{1}(R))$, and that $y$ is in the
essential image of $\Phi $.  Let $\alpha $ be a morphism from $\Phi x$
to $y$ for some point $(R,x)$ of $X_{0}$.  Consider the following
diagram.
\[
\begin{CD}
M_{\Phi x} @>\psi ^{M}_{\alpha }>> M_{y} @>\psi ^{M}_{\beta }>> M_{y'} \\
@V\tau _{x}VV @V\sigma _{y}VV @VV\sigma _{y'}V \\
N_{\Phi x} @>>\psi ^{N}_{\alpha }> N_{y} @>>\psi ^{N}_{\beta }> N_{y'}
\end{CD}
\]
By definition of $\sigma $, the left-hand square is commutative.  The
cocycle condition implies that $\psi _{\beta }\circ \psi _{\alpha }=\psi
_{\beta \alpha }$, so the definition of $\sigma $ also implies that the
outside square commutes.  Since the horizontal maps are isomorphisms,
the right-hand square must also be commutative.

We now extend the definition of $\sigma $ to an arbitrary point $(R,y)$
of $Y_{0}$.  The sheaf-theoretic essential image of $\Phi $ is all of
$Y_{0}$, we can choose a cover $\{R \xrightarrow{f_{i}} S_{i} \}$ of $R$
in the topology $\topology $ such that $y_{i}=Y_{0}(f_{i})(y)$ is in the
essential image of $\Phi $ for all $i$.  Let $y_{jk}$ denote the image
of $y$ in $Y_{0}(S_{j}\otimes _{R}S_{k})$.  We then have a commutative
diagram 
\[
\begin{CD}
M_{y} @>>> \prod M_{y_{i}} @>>> \prod M_{y_{jk}} \\
@. @V\prod \sigma _{y_{i}}VV @VV\prod \sigma _{y_{jk}} V \\
N_{y} @>>> \prod N_{y_{i}} @>>> \prod N_{y_{jk}}
\end{CD}
\]
where the right-hand horizontal maps are the difference of the two
restriction maps.  Thus each row expresses its left-hand entry as a
kernel.  The diagram commutes since $\sigma $ is a map of sheaves on the
essential image of $\Phi $.  Thus, there is a unique map $\sigma
_{y}\mathcolon M_{y}\xrightarrow{}N_{y}$ making the diagram commute.

We now check that $\sigma _{y}$ is independent of the choice of cover.
It suffices to show that $\sigma _{y}$ is unchanged if we replace the
cover $\{R\xrightarrow{}S_{i} \}$ by a refinement
$\{R\xrightarrow{}T_{j} \}$, since any two covers have a common
refinement.  If we denote the map coming from the refinement by
$\sigma '_{y}$, then we would have to have $\sigma '_{y_{i}}=\sigma
_{y_{i}}$, since some of the $T_{j}$ form a cover of $S_{i}$ and $\sigma
$ is a map of sheaves on the essential image of $\Phi $.  Then the sheaf
condition forces $\sigma' _{y}=\sigma _{y}$ as well.  In particular, if
$y$ is already in the essential image of $\Phi $, then we can take the
identity cover to find that the new definition of $\sigma $ is an
extension of our old definition.  

We now show that $\sigma $ is a map of sheaves over $Y_{0}$.  Suppose we
have a map $(R,y)\xrightarrow{f}(S,y')$ of points of $Y_{0}$.  Choose a
cover $\{R\xrightarrow{g_{i}}T_{i} \}$ of $R$ such that
$y_{i}=Y_{0}(g_{i})(y)$ is in the essential image of $\Phi $ for all
$i$.  Then there is an induced cover
$\{S\xrightarrow{h_{i}}U_{i}=S\otimes _{R}T_{i} \}$ of $S$.  The map $f$
induces corresponding maps $f_{i}\mathcolon (T_{i}, y_{i})
\xrightarrow{}(U_{i},y'_{i})$, where $y'_{i}=Y_{0}(h_{i})(y')$.  Since
$\sigma $ is a map of sheaves on the essential image of $\Phi $, we have
the commutative diagram below.
\[
\begin{CD}
M_{y_{i}} @>\sigma _{y_{i}}>> N_{y_{i}} \\
@VVV @VVV \\
M_{y'_{i}} @>>\sigma _{y'_{i}}>  N_{y'_{i}}
\end{CD}
\]
The sheaf condition and the definition of $\sigma $ then show that the
diagram
\[
\begin{CD}
M_{y} @>\sigma _{y}>> N_{y} \\
@VVV @VVV \\
M_{y'} @>>\sigma _{y'}> N_{y'}
\end{CD}
\]
is commutative, and so $\sigma $ is a map of sheaves over $Y_{0}$.  

The proof that $\sigma $ commutes with $\psi $, and so is a map of
sheaves over $(Y_{0},Y_{1})$, is similar.  
\end{proof}

Finally, we show that $\Phi ^{*}$ is essentially surjective.  

\begin{proposition}\label{prop-ess-surj}
Suppose $\Phi \mathcolon (X_{0},X_{1})\xrightarrow{}(Y_{0},Y_{1})$ is an
internal equivalence of presheaves of groupoids on $\site $.  Then $\Phi
^{*}\mathcolon \sheaves _{(Y_{0},Y_{1})}\xrightarrow{}\sheaves
_{(X_{0},X_{1})}$ is essentially surjective.
\end{proposition}

\begin{proof}
Suppose that $N$ is a sheaf over $(X_{0},X_{1})$.  We must construct a
sheaf $M$ over $(Y_{0},Y_{1})$ and an isomorphism $\Phi
^{*}M\xrightarrow{}N$ of sheaves.  We first construct $M_{y}$ for $y$ in
the essential image of $\Phi $, and show that it has the desired
properties there.  For every point $(R,y)$ in the
essential image of $\Phi $, choose a point $(R,x(y))$ of $X_{0}$ and a
morphism $\alpha (y)$ from $x(y)$ to $y$.  Note that this only requires
choosing over a set, since $\Aff $ is a small category.  Define
$M_{y}=N_{x(y)}$.  

We now construct the restriction of the structure map $\theta ^{M}$ to
the essential image of $\Phi $.  Suppose that we have a map
$(R,y)\xrightarrow{f}(S,y')$ between points of $Y_{0}$, where $(R,y)$ is
in the essential image of $\Phi $.  Let $\alpha '=Y_{1}(f)(\alpha (y))$,
so that $\alpha '$ is a morphism from $\Phi x'$ to $y'$, where
$x'=X_{0}(f)(x(y))$.  Then $\alpha (y')^{-1}\alpha '$ is a morphism from
$\Phi x'$ to $\Phi x(y')$.  Since $\Phi $ is full and faithful, there is
a unique morphism $\gamma $ of $X_{1}(S)$ from $x'$ to $x(y')$ such that
$\Phi \gamma =\alpha (y')^{-1}\alpha ' $, We then define $\theta
^{M}(f,y)\mathcolon M_{y}\xrightarrow{}M_{y'}$ to be the
composite
\[
M_{y}=N_{x(y)} \xrightarrow{\theta^{N} (f,x(y))}
N_{x'} \xrightarrow{\psi ^{N}_{\gamma }} N_{x(y')}=M_{y'}.
\]
We must check the functoriality conditions for $\theta ^{M}$ (restricted
to the essential image of $\Phi $).  First of
all, if $f$ is the identity map, then $\Phi \gamma $ will be the identity
morphism of $y$.  Since $\Phi $ is faithful, it follows that $\gamma $
is the identity morphism of $x(y)$.  The cocycle condition forces $\psi^{N}
_{\gamma }$ to be the identity map, and so $\theta ^{M}(1,y)$ is the
identity as required.  If $g\mathcolon (S,y')\xrightarrow{}(T,y'')$ is
another map of points of $Y_{0}$, a diagram chase involving the cocycle
condition for $\psi ^{N}$ and the fact that $\psi ^{N}$ is a map of
sheaves shows that $\theta ^{M} (gf,y)$ is the composition $\theta ^{M}
(g,y') \theta ^{M}(f,y)$.  

We now show that $M$ is a sheaf on the essential image of $\Phi
$.  Indeed, suppose $(R,y)$ is a point in the essential image of $\Phi
$, and $\{R\xrightarrow{}S_{i} \}$ is a cover of $R$ in $\topology $.
We must check that 
\[
M_{y}\xrightarrow{} \prod M_{y_{i}} \rightrightarrows \prod M_{y_{jk}}
\]
is an equalizer diagram.  We have an equalizer diagram 
\[
M_{y}=N_{x(y)} \xrightarrow{} \prod N_{x(y)_{i}} \rightrightarrows \prod
N_{x(y)_{jk}}
\]
since $N$ is a sheaf.  We construct an isomorphism from the bottom
diagram to the top, from which it follows that the top is also an
equalizer diagram.  The morphism $\alpha (y)\mathcolon \Phi
x(y)\xrightarrow{}y$ induces a morphism $\alpha (y)_{i}\mathcolon \Phi
x(y)_{i}\xrightarrow{}y_{i}$.  We also have the morphism $\alpha
(y_{i})\mathcolon \Phi x(y_{i})\xrightarrow{}y_{i}$.  The composition
$(\alpha (y_{i}))^{-1}\circ \alpha (y)_{i}=\Phi \gamma $ for a unique
$\gamma \mathcolon x(y)_{i}\xrightarrow{}x(y_{i})$, since $\Phi $ is
full and faithful.  Then $\psi _{\gamma }\mathcolon
N_{x(y)_{i}}\xrightarrow{}N_{x(y_{i})}=M_{y_{i}}$ defines the desired
isomorphism $\prod N_{x(y)_{i}}\xrightarrow{}\prod M_{y_{i}}$.  One
constructs the isomorphism $\prod N_{x(y)_{jk}}\xrightarrow{}\prod
N_{x(y_{jk})} = \prod M_{y_{jk}}$ in the same manner, using the
morphisms $\alpha (y)_{jk}\mathcolon \Phi x(y)_{jk}\xrightarrow{}y_{jk}$
and $\alpha (y_{jk})$.  The proof that the diagram below
\[
\begin{CD}
N_{x(y_{i})}=M_{y_{i}} @>>> N_{x(y_{ij})}=M_{y_{ij}} \\
@AAA @AAA \\
N_{x(y)_{i}} @>>> N_{x(y)_{ij}}
\end{CD}
\]
is commutative is a computation using the fact that $\psi ^{N}$ is a map
of sheaves, the cocycle condition, and the fact that $\Phi $ is
faithful.  

We now construct the restriction of the map $\psi ^{M}$ to the
essential image of $\Phi $.  Suppose $\beta $ is a morphism from $y$ to
$y'$, where $y$ is in the essential image of $\Phi $.  Then $\alpha
(y')^{-1}\beta \alpha (y)$ is a morphism from $\Phi x(y)$ to $\Phi
x(y')$.  Since $\Phi $ is full and faithful, there is a unique morphism
$\gamma $ from $x(y)$ to $x(y')$ such that $\Phi \gamma =\alpha
(y')^{-1}\beta \alpha (y)$.  Hence we can define $\psi ^{M}_{\beta
}=\psi ^{N}_{\gamma }$.  We leave to the reader the diagram chase
showing that $\psi $ is a map of sheaves.  

We now construct the desired isomorphism of sheaves $\tau \mathcolon
\Phi ^{*}M \xrightarrow{}N$.  (Since $\Phi ^{*}M$ is determined by the
restriction of $M$ to the image of $\Phi $, we can do this even though
we have not completed the definition of $M$).  Suppose $(R,x)$ is a
point of $X_{0}$.  Then $\alpha (\Phi x)$ is a morphism from $\Phi
(x(\Phi x))$ to $\Phi x$.  Since $\Phi $ is full and faithful, there is
a unique morphism $\beta $ from $x(\Phi x)$ to $x$ such that $\Phi \beta
=\alpha (\Phi x)$.  We define 
\[
\tau _{x}=\psi ^{N}_{\beta }\mathcolon M_{\Phi x}=N_{x(\Phi x)}
\xrightarrow{} N_{x}.  
\]
Obviously $\tau _{x}$ is an isomorphism, but we must check it is
compatible with the structure maps.  We leave these checks to the
reader; both are diagram chases.  

We have now defined a sheaf $M$ on the essential image of $\Phi $, and
to complete the proof we need only extend it to a sheaf on all of
$(Y_{0},Y_{1})$.  For each point $(R,y)$ of $Y_{0}$, choose a cover
$C(y)=\{R\xrightarrow{f_{i}}S_{i} \}$ such that $y_{i}=Y_{0}(f_{i})(y)$
is in the essential image of $\Phi $ for all $i$, making sure to choose
the identity cover when $y$ is already in the essential image of $\Phi
$.  Once again, we can do this since $\Aff $ is a small category.  We
then define $M_{y}$ as we must if we are going to get a sheaf, as the
equalizer of the two maps of $R$-modules
\[
\prod _{i} M_{y_{i}} \rightrightarrows \prod _{jk} M_{y_{jk}}.
\]

This definition of $M_{y}$ will of course depend on the choice of
cover $C(y)$.  Suppose $D=\{R\xrightarrow{}T_{m} \}$ is some other
cover such that $y_{m}$ is in the essential image of $\Phi $ for all
$m$.   We claim that there is a canonical equalizer diagram 
\[
M_{y} \xrightarrow{} \prod M_{y_{m}} \rightrightarrows \prod M_{y_{np}}.
\]
To see this, let $M^{D}_{y}$ denote the pullback of the two arrows 
\[
\prod _{m} M_{y_{m}} \rightrightarrows \prod _{np} M_{y_{np}}.  
\]
We claim that there is a canonical isomorphism
$M^{D}_{y}\xrightarrow{}M_{y}$.  It suffices to check this when $D$ is a
refinement of $C(y)$, since any two covers have a common refinement.  In
this case, there is a diagram 
\[
M_{y}\xrightarrow{} \prod _{m} M_{y_{m}} \rightrightarrows \prod _{np}
M_{y_{np}},
\]
where the first map is induced by first mapping to $M_{y_{i}}$, and then
using the structure maps of $M$ restricted to the essential image of
$\Phi $ to map further to $M_{y_{m}}$.  It suffices to prove that this
diagram is an equalizer.  It is easy to check that $M_{y}$ maps into
the equalizer.  If $t\in M_{y}$ maps to $0$ in each $M_{y_{m}}$, then,
using the fact that $M$ restricted to the essential image of $\Phi $ is
a sheaf, we find that $t$ maps to $0$ in each $M_{y_{i}}$.  By
definition of $M_{y}$, then, $t=0$.  Similarly, suppose $(t_{m})\in \prod
M_{y_{m}}$ is in the equalizer.  Again using the fact that $M$
restricted to the essential image of $\Phi $ is a sheaf, we construct
an element $(t_{i})\in \prod M_{y_{i}}$.  The images
of $t_{i}$ and $t_{j}$ in $M_{y_{ij}}$ coincide, since they coincide
after restriction to the induced cover.  Thus we get an element $t\in
M_{y}$ restricting to the $t_{i}$.  It follows that $t$ restricts to the
$t_{m}$ as well, and so $M_{y}$ is the desired equalizer.  

Now we can construct the structure maps of $M$.  Suppose
$(R,y)\xrightarrow{}(S,z)$ is a map of points of $Y_{0}$.  The cover
$C(y)=\{R\xrightarrow{}S_{i} \}$ of $R$ induces a cover
$D=\{S\xrightarrow{}S\otimes _{R}S_{i} \}$ of $S$, and the restriction
$z_{i}$ of $z$ is in the essential image of $\Phi $ for all $i$, since
$y_{i}$ is so.  Thus we get a map from 
\[
\prod M_{y_{i}} \rightrightarrows \prod M_{y_{jk}} 
\]
to 
\[
\prod M_{z_{i}} \rightrightarrows \prod M_{z_{jk}},
\]
and so an induced map $M_{y}\xrightarrow{}M^{D}_{z}$ on the equalizers.
After composing this with the canonical isomorphism
$M^{D}_{z}\xrightarrow{}M_{z}$, we get the desired structure map $\theta
\mathcolon M_{y}\xrightarrow{}M_{z}$.  Since we chose the identity cover
when $y$ was already in the essential image of $\Phi $, this extends the
definition we have already given in that case.  We leave it to the
reader to check the functoriality of $\theta $.  

We now show that $M$ is a sheaf.  Suppose $(R,y)$ is a point of $Y_{0}$
and $\{(R,y)\xrightarrow{}(T_{m},y_{m}) \}$ is a cover of $R$.  Let
$C(y)=\{(R,y)\xrightarrow{}(S_{i},y_{i}) \}$ be the given cover of $R$,
so that each $y_{i}$ is in the essential image of $\Phi $.  Then
$\{S_{i}\xrightarrow{}T_{m}\otimes _{R}S_{i} \}$ is a cover of $S_{i}$,
and each $y_{mi}$ is the essential image of $\Phi $ since each $y_{i}$
is.  Similarly, $\{T_{m}\xrightarrow{}T_{m}\otimes _{R}S_{i} \}$ is a
cover of $T_{m}$.  Thus we get the commutative diagram below.  
\[
\begin{CD}
M_{y} @>>> \prod _{m} M_{y_{m}} @>>> \prod _{np} M_{y_{np}} \\
@VVV @VVV @VVV \\
\prod _{i} M_{y_{i}} @>>> \prod _{mi} M_{y_{mi}} @>>> \prod _{npi}
M_{y_{npi}} \\
@VVV @VVV @VVV \\
\prod _{jk} M_{y_{jk}} @>>> \prod _{mjk} M_{y_{mjk}} @>>> \prod _{npjk}
M_{y_{npjk}} 
\end{CD}
\]
The subscripts $m$, $n$, and $p$ all refer to the $T_{m}$, and the
subscripts $i,j$ and $k$ all refer to the $S_{i}$.  So, for example,
$y_{npi}$ is the image of $y$ in $Y_{0}(T_{n}\otimes _{R}T_{p}\otimes
_{R}S_{i})$.  The right-hand horizontal arrows are all the differences
of the two restriction maps.  This means that the second and third rows
express their left-hand entries as kernels, since $M$ restricted to the
essential image of $\Phi $ is a sheaf.  Similarly, the bottom vertical
arrows are also differences of the two restriction maps.  It follows
that each column expresses its top entry as a kernel, since the
definition of $M$ does not depend on which cover we choose, up to
isomorphism.  A diagram chase then shows that the top row expresses
$M_{y}$ as a kernel, which means that $M$ is a sheaf.  

We now construct the isomorphism $\psi \mathcolon \dom
^{*}M\xrightarrow{}\codom ^{*}M$.  Suppose $\alpha \mathcolon
y\xrightarrow{}z$ is a morphism in $Y_{1}(R)$.  Let
$\{R\xrightarrow{}S_{i} \}$ be the given cover of $(R,y)$, so that each
$y_{i}$ is in the essential image of $\Phi $.  It follows that $z_{i}$
is also in the essential image of $\Phi $ for all $i$.  Let $\alpha
_{i}\mathcolon y_{i}\xrightarrow{}z_{i}$ denote the image of $\alpha $
in $Y_{1}(S_{i})$, and similarly let $\alpha _{jk}$ denote the image of
$\alpha $ in $Y_{1}(S_{j}\otimes _{R}S_{k})$.  Then we have a
commutative diagram
\[
\begin{CD}
M_{y} @>>> \prod M_{y_{i}} @>>> \prod M_{y_{jk}} \\
@.  @V\prod \psi_{\alpha _{i}} VV @VV\prod \psi _{\alpha _{jk}} V \\
M_{z} @>>> \prod M_{z_{i}} @>>> \prod M_{z_{jk}}.
\end{CD}
\]
Here the right-hand horizontal arrows are differences of restriction
maps, as usual.  The top row is an equalizer by definition, and we have
proved that the bottom row is also an equalizer diagram.  Hence there is
a unique map $\psi _{\alpha }\mathcolon M_{y}\xrightarrow{}M_{z}$,
necessarily an isomorphism, making the diagram commute.  The facts that
$\psi $ satisfies the cocycle condition and is a map of sheaves are the
usual sheaf-theoretic diagram chases, and we leave them to the reader.  
\end{proof}

\section{Quasi-coherent sheaves}\label{sec-quasi}

The object of this section is to prove Theorem~\ref{main-C}, showing
that if $\Phi \mathcolon (X_{0},X_{1})\xrightarrow{}(Y_{0},Y_{1})$ is an
internal equivalence of presheaves of groupoids in the flat topology,
then $\Phi ^{*}\mathcolon \qcsheaves _{(Y_{0},Y_{1})}\xrightarrow{}
\qcsheaves _{(X_{0},X_{1})}$ is an equivalence of categories of
quasi-coherent sheaves.  This theorem can be viewed as a manifestation
of faithfully flat descent; we have seen already that $\Phi
^{*}\mathcolon \sheaves _{(Y_{0},Y_{1})}\xrightarrow{}\sheaves
_{(X_{0},X_{1})}$ is an equivalence of categories, and we use faithfully
flat descent to conclude that quasi-coherent sheaves are a full
subcategory of sheaves in the flat topology.

Recall that a cover of $R$ in the flat, or fpqc, topology is a
\emph{finite} collection of maps $\{R\xrightarrow{}S_{i} \}$ such that
each $S_{i}$ is flat over $R$, and the product $\prod S_{i}$ is
faithfully flat over $R$.  This also defines the flat topology on $\Aff
_{*}$.  

We use faithfully flat descent in the form of the following well-known
lemma. 

\begin{lemma}\label{lem-descent}
Suppose $\{R\xrightarrow{}S_{i} \}$ is a cover of $R$ in the flat
topology on $\Aff $, and $M$ is an $R$-module.  Then the diagram
\[
M \xrightarrow{} \prod _{i} S_{i}\otimes _{R} M \rightrightarrows \prod
_{jk} S_{j}\otimes _{R} S_{k} \otimes _{R} M
\]
is an equalizer in the category of $R$-modules.  
\end{lemma}

Of course, the two maps in the equalizer take $s\otimes m\in
S_{i}\otimes M$ to $(1\otimes s_{i}\otimes m)\in \prod _{ji}
S_{j}\otimes _{R}S_{i}\otimes _{R}M$ and to $s_{i}\otimes 1\otimes m\in
\prod _{ik}S_{i}\otimes _{R}S_{k}\otimes _{R}M$.  

As usual, this lemma also works in the graded case, with the same
proof. 

\begin{proof}
Let $S=\prod_{i} S_{i}$.  Since the product is finite, it suffices to
show that 
\[
M \xrightarrow{} S\otimes _{R} M \rightrightarrows S\otimes _{R} S
\otimes _{R} M
\]
is an equalizer for all $R$-modules $M$.  Since $S$ is faithfully flat,
it suffices to show that 
\[
S\otimes _{R} M \xrightarrow{} S\otimes _{R}S\otimes _{R} M
\rightrightarrows S \otimes _{R} S\otimes _{R} S\otimes _{R} M
\]
is an equalizer for all $M$.  But, before tensoring with $M$, this
sequence is just the beginning of the bar resolution of $S$ as an
$R$-algebra; since the bar resolution is contractible, this diagram
remains an equalizer after tensoring with $M$.  
\end{proof}

Lemma~\ref{lem-descent} leads immediately to the following proposition,
which is also true in the graded case. 

\begin{proposition}\label{prop-descent}
Suppose $M$ is a quasi-coherent sheaf over a presheaf of groupoids
$(X_{0},X_{1})$ on $\Aff $.  Then $M$ is a sheaf in the flat topology.  
\end{proposition}

\begin{proof}
Suppose $(R,y)$ is a point of $X_{0}$, and
$\{(R,y)\xrightarrow{}(S_{i},y_{i}) \}$ is a cover in the flat
topology.  We must show that the diagram 
\[
E_{y}= (M_{y}\xrightarrow{}\prod M_{y_{i}}\rightrightarrows \prod M_{y_{jk}})
\]
is an equalizer diagram.  But, since $M$ is quasi-coherent, $E_{y}$ is
isomorphic to the diagram 
\[
M_{y} \xrightarrow{} \prod S_{i}\otimes _{R} M_{y} \rightrightarrows
\prod S_{j}\otimes _{R}S_{k}\otimes _{R}M_{y},
\]
which is an equalizer diagram by Lemma~\ref{lem-descent}.  
\end{proof}

We will also need a lemma about purity of equalizer diagrams.  

\begin{definition}\label{defn-pure}
Suppose $E$ is an equalizer diagram of the form 
\[
A \xrightarrow{} B \rightrightarrows C
\]
in the category of $R$-modules for some commutative ring $R$.  We say
that $E$ is \emph{pure} if $S\otimes _{R} E$ is still an equalizer
diagram for all commutative $R$-algebras $S$.  
\end{definition}

One can also define purity using arbitrary $R$-modules $S$.  We prefer
this definition because it is the concept we need, but in fact the two
definitions are equivalent.  Either definition also works in the graded
case with the obvious changes.  

\begin{lemma}\label{lem-purity}
Suppose $E$ is an equalizer diagram of $R$-modules for some commutative
ring $R$.  Suppose $\{S_{i} \}$ is a set of flat commutative $R$-algebra
such that $S_{i}\otimes _{R} E$ is pure for all $i$ and $S=\bigoplus
_{i}S_{i}$ is faithfully flat over $R$.  Then $E$ is pure.  
\end{lemma}

\begin{proof}
Suppose $T$ is an arbitrary $R$-algebra.  Then $(T\otimes
_{R}S_{i})\otimes _{S_{i}}(S_{i}\otimes _{R}E)$ is an equalizer diagram
since $S_{i}\otimes _{R}E$ is pure, but
\[
(T\otimes _{R}S_{i})\otimes _{S_{i}}(S_{i}\otimes _{R}E) \cong (T\otimes
_{R}S_{i}) \otimes _{T}(T\otimes _{R}E).
\]
Thus $(T\otimes _{R}S)\otimes _{T}(T\otimes _{R}E)$ is also an equalizer
diagram, being a direct sum of equalizer diagrams.  Since $T\otimes
_{R}S$ is faithfully flat over $T$, it follows that $T\otimes _{R}E$ is
an equalizer diagram.
\end{proof}

We can now prove that quasi-coherent sheaves are homotopy invariant in
the flat topology.  The following theorem is Theorem~\ref{main-C} of the
introduction.  

\begin{theorem}\label{thm-qc-ess-surj}
Suppose $\Phi \mathcolon (X_{0},X_{1})\xrightarrow{}(Y_{0},Y_{1})$ is an
internal equivalence of presheaves of groupoids on $\site $, where
$\topology $ is the flat topology.  Then $\Phi ^{*}\mathcolon \qcsheaves
_{(Y_{0},Y_{1})}\xrightarrow{}\qcsheaves _{(X_{0},X_{1})}$ is an
equivalence of categories.
\end{theorem}

This theorem is also true in the graded case, with the same proof.  

\begin{proof}
Since $\Phi ^{*}\mathcolon \sheaves
_{(Y_{0},Y_{1})}\xrightarrow{}\sheaves _{(X_{0},X_{1})}$ is an
equivalence of categories, and quasi-coherent sheaves are a full
subcategory of sheaves in the flat topology by
Proposition~\ref{prop-descent}, we find immediately that $\Phi
^{*}\mathcolon \qcsheaves _{(Y_{0},Y_{1})}\xrightarrow{}\qcsheaves
_{(X_{0},X_{1})}$ is full and faithful.  It remains to show that it is
essentially surjective. 

Suppose $N$ is a quasi-coherent sheaf over $(X_{0},X_{1})$.  Because
$\Phi ^{*}\mathcolon \sheaves _{Y_{0},Y_{1}}\xrightarrow{}\sheaves
_{(X_{0},X_{1})}$ is an equivalence of categories, there is a sheaf $M$
in the flat topology, over $(Y_{0},Y_{1})$, such that $\Phi ^{*}M\cong
N$.  We will show that $M$ is in fact quasi-coherent, so that $\Phi
^{*}$ is essentially surjective on quasi-coherent sheaves.  To do so, we
must show that, if $(R,y)\xrightarrow{f}(S,y')$ is a map of points of
$Y_{0}$, then the adjoint $S\otimes _{R}M_{y}\xrightarrow{\rho
^{M}(f)}M_{y'}$ of the structure map of $M$ is an isomorphism.

First suppose that $y$ is in the essential image of $\Phi $.  Then there
is an $x\in X_{0}(R)$ and a map $\alpha \mathcolon \Phi
x\xrightarrow{}y$.  Let $x'=f(x)\in X_{0}(S)$, so that $f(\alpha
)=X_{1}(f)(\alpha )\mathcolon \Phi x'\xrightarrow{}z$.  Then we have the
commutative diagram below.
\[
\begin{CD}
S\otimes _{R} N_{x} @>\rho ^{N}(f)>> N_{x'} \\
@V\cong VV @VV\cong V \\
S\otimes _{R} M_{\Phi x} @>\rho ^{M}(f)>> M_{\Phi x'}  \\
@V1\otimes \psi _{\alpha }VV @VV\psi _{f\alpha } V \\
S\otimes _{R}M_{y} @>>\rho ^{M}(f)> M_{y'}
\end{CD}
\]
The top square of this diagram commutes because $\Phi ^{*}M\cong N$ as
sheaves, and the bottom square commutes because $\psi $ is a map of
sheaves.  The vertical maps are isomorphisms, and the top horizontal map
is an isomorphism since $N$ is quasi-coherent.  Hence the bottom
horizontal map is an isomorphism as well.  

In fact, if $y$ is in the essential image of $\Phi $ and
$\{R\xrightarrow{}S_{i} \}$ is a cover of $R$ in the flat topology, we
claim that the equalizer diagram 
\begin{equation}\label{eq-eq}
E= E_{y}=(M_{y} \xrightarrow{} \prod M_{y_{i}} \rightrightarrows \prod
M_{y_{jk}})
\end{equation}
is pure.  Indeed, suppose $S$ is an $R$-algebra, so we have $f\mathcolon
(R,y)\xrightarrow{}(S,y')$.  Then $\{S\xrightarrow{} S\otimes
_{R}S_{i}\}$ is a cover of $S$ in the flat topology.  It follows from
what we have just done (and the fact that covers in the flat topology
are finite), that the diagram $S\otimes _{R}E_{y}$ is isomorphic to
$E_{y'}$, and so is still an equalizer diagram.

Now suppose $y$ is an arbitrary point of $Y_{0}$.  Since the
sheaf-theoretic essential image of $\Phi $ is all of $Y_{0}$, we can
choose a cover $\{R\xrightarrow{}S_{i} \}$ such that each $y_{i}$ is in
the essential image of $\Phi $.  There is an induced cover
$\{S\xrightarrow{}S\otimes _{R}S_{i} \}$ of $S$, and maps
$f_{i}\mathcolon (S_{i},y_{i})\xrightarrow{}(S\otimes
_{R}S_{i},y'_{i})$, so each $y'_{i}$ is also in the essential image of
$\Phi $.  We then get the commutative diagram below, which is a map from
the diagram $S\otimes _{R}E_{y}$ to $E_{z}$.  
\[
\begin{CD}
S\otimes _{R} M_{y} @>>>\prod S\otimes _{R}M_{y_{i}} @>1\otimes d>>
\prod S\otimes 
_{R}M_{y_{jk}} \\
@V\rho _{f} VV @V\prod \rho (f_{i})VV @VV\prod \rho (f_{jk})V \\
M_{z} @>>> \prod M_{z_{i}} @>>d> \prod M_{z_{jk}}.
\end{CD}
\]
Here the map $d$ is the difference between the two restriction maps, so
the bottom row expresses $M_{z}$ as a kernel.  We have already seen that
the maps $\rho (f_{i})$ and $\rho (f_{jk})$ are isomorphisms, so if we
knew that $S\otimes _{R}E_{y}$ were an equalizer diagram, we would be
able to conclude that $\rho (f)$ is an isomorphism, and therefore that
$M$ is quasi-coherent.

In particular, if $S$ is flat over $R$, we conclude that the diagram
$S\otimes _{R}E_{y}$ is isomorphic to the equalizer diagram $E_{y'}$.
In case $y'$ is in the essential image of $\Phi $, we have proved that
$E_{y'}$ is pure.  In particular, $S_{i}\otimes
_{R}E$ is a pure equalizer diagram for all $i$.  Since $\prod S_{i}$ is
faithfully flat over $R$, it follows from Lemma~\ref{lem-purity} that
the equalizer diagram $E$ is pure.  Thus, for any $S$, $S\otimes _{R}E$
is an equalizer diagram, and so $M$ is quasi-coherent.  
\end{proof}

\section{Hopf algebroids}\label{sec-Hopf}

In this section, we prove Theorem~\ref{main-D} of the introduction,
characterizing those maps of Hopf algebroids which induce internal
equivalences in the flat topology of the corresponding presheaves of
groupoids.  

Suppose $f=(f_{0},f_{1})\mathcolon (A,\Gamma )\xrightarrow{}(B,\Sigma )$
is a map of Hopf algebroids.  See~\cite[Definition~A1.1.7]{ravenel} for
an explicit definition of this, though of course $f$ is 
equivalent to a map $\Phi =f^{*}\mathcolon (\Spec B,\Spec \Sigma
)\xrightarrow{}(\Spec A,\Spec \Gamma )$ of sheaves of groupoids on $\Aff
$.  A map of Hopf algebroids induces a map 
\[
B\otimes _{A} \Gamma \otimes _{A} B \xrightarrow{\eta _{L}\otimes
f_{1}\otimes \eta _{R}} \Sigma _{\eta _{L}f_{0}} \otimes _{A}
\bimod{\Sigma }{f_{1}\eta _{L}}{f_{1}\eta _{R}} \otimes _{A}
\bimod{\Sigma }{\eta _{R}f_{0}}{\ } \xrightarrow{\mu } \Sigma ,
\]
where $\mu $ denotes multiplication.  Note that $\mu $ makes sense since
$f_{1}\eta _{L}=\eta _{L}f_{0}$ and $f_{1}\eta _{R}=\eta _{R}f_{0}$.  By
abuse of notation, we denote this map simply by $\eta _{L}\otimes
f_{1}\otimes \eta _{R}$.  

Our goal is to characterize those $f$ for which $f^{*}$ is a weak
equivalence.  We begin by determining when $f^{*}$ is faithful.  

\begin{proposition}\label{prop-Hopf-faithful}
Suppose $f=(f_{0},f_{1})\mathcolon (A,\Gamma )\xrightarrow{}(B,\Sigma )$
is a map of Hopf algebroids.  Then $f^{*}\mathcolon (\Spec B,\Spec
\Sigma )\xrightarrow{}(\Spec A,\Spec \Gamma )$ is faithful if and only
if $\eta _{L}\otimes f_{1}\otimes \eta _{R}\mathcolon B\otimes _{A}
\Gamma \otimes _{A} B\xrightarrow{}\Sigma $ is an epimorphism in $\Ring
$.  
\end{proposition}

Recall that an epimorphism in $\Ring $ need not be surjective; the map
from the integers to the rational numbers is a ring epimorphism.  Also
note that the obvious generalization of this proposition holds for
graded Hopf algebroids.  

\begin{proof}
Given $\alpha, \beta  \mathcolon \Sigma \xrightarrow{}R$, 
\[
\alpha \circ
(\eta _{L}\otimes f_{1}\otimes \eta _{R})=\beta \circ (\eta _{L}\otimes
f_{1}\otimes \eta _{R})
\]
if and only if $\alpha $ and $\beta $ have the same domain and codomain
when thought of as morphisms of $(\Spec B,\Spec \Sigma )(R)$ and 
$f^{*}\alpha =f^{*}\beta $.  The proposition follows.  
\end{proof}

We now determine when $f^{*}$ is full.  

\begin{proposition}\label{prop-Hopf-full}
Suppose $f=(f_{0},f_{1})\mathcolon (A,\Gamma )\xrightarrow{}(B,\Sigma )$
is a map of Hopf algebroids.  Then $f^{*}\mathcolon (\Spec B,\Spec
\Sigma )\xrightarrow{}(\Spec A,\Spec \Gamma )$ is full if and only
if $\eta _{L}\otimes f_{1}\otimes \eta _{R}\mathcolon B\otimes _{A}
\Gamma \otimes _{A} B\xrightarrow{}\Sigma $ is a split monomorphism of
rings.  
\end{proposition}

Once again, the obvious generalization of this proposition is true in
the graded case.  

\begin{proof}
The map $f^{*}$ is full if and only if every morphism 
\[
\beta \mathcolon f^{*}x\xrightarrow{}f^{*}y \in (\Spec A,\Spec \Gamma
)(R)
\]
is equal to $f^{*}\alpha $ for some morphism $\alpha \mathcolon
x\xrightarrow{}y$ of $(\Spec B,\Spec \Sigma )(R)$.  Said another way,
$f^{*}$ is full if and only if every ring homomorphism
\[
x\otimes \beta \otimes y\mathcolon B\otimes _{A}\Gamma \otimes
_{A}B\xrightarrow{}R 
\]
can be extended through $\eta _{L}\otimes f_{1}\otimes \eta _{R}$ to a
ring homomorphism $\Sigma \xrightarrow{}R$.  This is equivalent to $\eta
_{L}\otimes f_{1}\otimes \eta _{R}$ being a split monomorphism.  
\end{proof}

\begin{corollary}\label{cor-Hopf-full}
Suppose $f=(f_{0},f_{1})\mathcolon (A,\Gamma )\xrightarrow{}(B,\Sigma )$
is a map of Hopf algebroids.  Then $f^{*}\mathcolon (\Spec B,\Spec
\Sigma )\xrightarrow{}(\Spec A,\Spec \Gamma )$ is fully faithful if and
only if $\eta _{L}\otimes f_{1}\otimes \eta _{R}\mathcolon B\otimes _{A}
\Gamma \otimes _{A} B\xrightarrow{}\Sigma $ is an isomorphism.  
\end{corollary}

\begin{proof}
Any map $g\mathcolon R\xrightarrow{}S$ of rings that is both a split
monomorphism and a ring epimorphism is an isomorphism.  Indeed, $\Ring
(g,T)\mathcolon \Ring (S,T)\xrightarrow{}\Ring (R,T)$ is monic since $g$
is a ring epimorphism and epic since $g$ is a split monomorphism, so is
an isomorphism for all $T$.  
\end{proof}

Finally, we need to determine the sheaf-theoretic essential image of
$f^{*}$ is all of $\Spec A$.  For this we need the map $f_{0}\otimes \eta
_{R}\mathcolon A\xrightarrow{} B \otimes _{A} \Gamma $ defined as the
composite
\[
A \cong A\otimes _{A} A \xrightarrow{f_{0}\otimes \eta _{R}} B\otimes _{A}
\Gamma . 
\]

\begin{proposition}\label{prop-Hopf-essentially}
Suppose $f=(f_{0},f_{1})\mathcolon (A,\Gamma )\xrightarrow{}(B,\Sigma )$
is a map of Hopf algebroids.  Then the sheaf-theoretic essential image
of 
\[
f^{*}\mathcolon (\Spec B,\Spec \Sigma )\xrightarrow{}(\Spec A,\Spec
\Gamma )
\]
is all of $\Spec A$ if and only if there is a ring map $g\mathcolon
B\otimes _{A}\Gamma \xrightarrow{}C$ such that $g(f_{0}\otimes \eta
_{R})$ exhibits $C$ as a faithfully flat extension of $A$.
\end{proposition}

This proposition is also true in the graded case, with the same proof. 

\begin{proof}
We first determine when $y\mathcolon A\xrightarrow{}R$ is in the
essential image of $f^{*}$.  For this to happen we need an object
$x\mathcolon B\xrightarrow{}R$ and a morphism $\alpha \mathcolon \Gamma
\xrightarrow{}R$ from $f^{*}x$ to $y$.  A morphism $\alpha $ from
$f^{*}x$ to anywhere is equivalent to the composite
\[
B\otimes _{A}\Gamma \xrightarrow{x\otimes \alpha } R_{xf_{0}} \otimes
_{A} \bimod{R}{\alpha \eta _{L}}{\ } \xrightarrow{\mu } R,
\]
which we also denote, by abuse of notation, by $x\otimes \alpha $.  The
codomain of $\alpha $ is the composite $(x\otimes \alpha )(f_{0}\otimes
\eta _{R})\mathcolon A\xrightarrow{}R$.  Altogether then, $y$ is in the
essential image of $f^{*}$ if and only if there is a map $h\mathcolon
B\otimes _{A}\Gamma $ such that $h(f_{0}\otimes \eta _{R})=y$.

Now, suppose the sheaf-theoretic essential image of $f^{*}$ is all of
$\Spec A$.  Then there must be a cover $\{A\xrightarrow{h_{i}}S_{i} \}$
such that the image of the identity map of $A$, namely $h_{i}$, is in
the essential image of $f^{*}$ for all $i$.  By the preceding paragraph,
this is true if and only if there exist maps $g_{i}\mathcolon B\otimes
_{A} \Gamma \xrightarrow{}S_{i}$ such that $g_{i}(f_{0}\otimes \eta
_{R})=h_{i}$.  Let $C$ be the product of the $S_{i}$ and let
$g\mathcolon B\otimes _{A}\Gamma \xrightarrow{}C$ be the product of the
$g_{i}$.  Then $g(f_{0}\otimes \eta _{R})$ is the product of the
$h_{i}$, which displays $C$ as a faithfully flat extension of $A$ since
$\{A\xrightarrow{h_{i}}S_{i} \}$ is a cover of $A$.

Conversely, suppose there is a ring map $g\mathcolon B\otimes _{A}\Gamma
\xrightarrow{}C$ such that $h=g(f_{0}\otimes \eta _{R})$ exhibits $C$ as
a faithfully flat extension of $A$.  Suppose $y\mathcolon
A\xrightarrow{}R$ is an arbitrary point of $(\Spec A,\Spec \Gamma )(R)$.
Then
\[
R\cong A\otimes _{A}R\xrightarrow{h\otimes 1}C\otimes _{A}R
\]
is a cover of $R$.  One can easily check that the image of $y$ in
$(\Spec A,\Spec \Gamma )(C\otimes _{A}R)$ is the composite
\[
A\xrightarrow{h} C \cong C\otimes _{A}A \xrightarrow{1\otimes y}
C\otimes _{A}R.
\]
Since $h=g(f_{0}\otimes \eta _{R})$, the image of $y$ is in the essential
image of $f^{*}$, and so $y$ is in the sheaf-theoretic essential image
of $f^{*}$.
\end{proof}

Note that the proof of Proposition~\ref{prop-Hopf-essentially} can be
easily modified to prove the known result that $f^{*}$ is essentially
surjective if and only if $f_{0}\otimes \eta _{R}\mathcolon
A\xrightarrow{}B\otimes _{A}\Gamma $ is a split monomorphism.  

Altogether then, we have the following theorem, which is
Theorem~\ref{main-D} of the introduction. 

\begin{theorem}\label{thm-Hopf}
Suppose $f=(f_{0},f_{1})\mathcolon (A,\Gamma )\xrightarrow{}(B,\Sigma )$
is a map of Hopf algebroids.  Then $f^{*}\mathcolon (\Spec B,\Spec
\Sigma )\xrightarrow{}(\Spec A,\Spec \Gamma )$ is an internal
equivalence in the flat topology if and only if
\[
\eta _{L}\otimes f_{1}\otimes \eta _{R}\mathcolon B\otimes _{A}\Gamma
\otimes _{A}B \xrightarrow{} \Sigma 
\]
is an isomorphism and there is a ring map $g\mathcolon B\otimes
_{A}\Gamma \xrightarrow{}C$ such that $g(f_{0}\otimes \eta _{R})$
exhibits $C$ as a faithfully flat extension of $A$.
\end{theorem}

This characterization of internal equivalences shows in particular that
$\Sigma $ is determined by $(A,\Gamma )$ and $f_{0}$.  In fact, if $(A,
\Gamma )$ is any Hopf algebroid, and $f\mathcolon A\xrightarrow{}B$ is a
ring homomorphism, there is a unique (up to isomorphism) Hopf algebroid
$(B, \Gamma _{f})$ and map of Hopf algebroids $(f,f_{1})$ such that the
map $\eta _{L}\otimes f_{1}\otimes \eta _{R}$ is an isomorphism.  To
show existence, we take $\Gamma _{f}=B\otimes _{A}\Gamma \otimes _{A}B$
and define the structure maps as follows:
\begin{align*}
\eta _{L} & \mathcolon B \cong B\otimes _{A}A\otimes _{A}A
\xrightarrow{1\otimes \eta _{L}\otimes f} B\otimes _{A} \Gamma \otimes
_{A} B\usc \\
\eta _{R} & \mathcolon B \cong A \otimes _{A} A \otimes _{A} B
\xrightarrow{f\otimes \eta _{R} \otimes 1} B\otimes _{A} \Gamma \otimes
_{A} B\usc \\
\epsilon & \mathcolon B\otimes _{A} \Gamma \otimes _{A} B
\xrightarrow{1\otimes \epsilon \otimes 1} B \otimes _{A} A \otimes _{A}
B \cong B \otimes _{A} B \xrightarrow{\mu } B\usc \\
c & \mathcolon B\otimes _{A} \Gamma \otimes _{A} B \xrightarrow{1\otimes
c\otimes 1} B \otimes _{A} \bimod{\Gamma }{\eta _{R}}{\eta _{L}} \otimes
_{A} B \xrightarrow{\tau } B \otimes _{A} \bimod{\Gamma }{\eta
_{L}}{\eta _{R}} \otimes _{A} B\usc \\
\Delta & \mathcolon B \otimes _{A} \Gamma \otimes _{A} B
\xrightarrow{1\otimes \Delta \otimes 1} B \otimes _{A} \Gamma \otimes
_{A} \Gamma \otimes _{A} B \cong B \otimes _{A} \Gamma \otimes _{A} A
\otimes _{A} \Gamma \otimes _{A} B \\
& \xrightarrow{1\otimes 1\otimes
f\otimes 1\otimes 1} B \otimes _{A} \Gamma \otimes _{A} B \otimes _{A}
\Gamma \otimes _{A} B \cong (B\otimes _{A}\Gamma \otimes _{A}B)\otimes
_{B}(B\otimes _{A}\Gamma \otimes _{A}B).
\end{align*}
We leave it to the reader to check that this does define a Hopf
algebroid.  We define $f_{1}\mathcolon \Gamma \xrightarrow{}\Gamma _{f}$
to be the composite 
\[
\Gamma \cong A\otimes _{A}\Gamma \otimes _{A} A \xrightarrow{f\otimes
1\otimes f} B \otimes _{A} \Gamma \otimes _{A} B.  
\]
We leave it to the reader to check that this defines a map of Hopf
algebroids, and also to check our uniqueness claims.  

We therefore have the following corollary. 

\begin{corollary}\label{cor-Hopf-weak}
Suppose $f=(f_{0},f_{1})\mathcolon (A,\Gamma )\xrightarrow{}(B,\Sigma )$
is a map of Hopf algebroids.  Then $f^{*}\mathcolon (\Spec B,\Spec
\Sigma )\xrightarrow{}(\Spec A,\Spec \Gamma )$ is an internal
equivalence in the flat topology if and only if $(B,\Sigma )$ is
isomorphic over $(A,\Gamma )$ to $(B,\Gamma _{f_{0}})$ and there is a
ring map $g\mathcolon B\otimes _{A}\Gamma \xrightarrow{}C$ such that
$g(f_{0}\otimes \eta _{R})$ exhibits $C$ as a faithfully flat extension
of $A$.
\end{corollary}

The conditions in Corollary~\ref{cor-Hopf-weak} have appeared before,
in~\cite[Theorem~3.3]{hovey-sadofsky-picard} and in~\cite{hopkins-hopf}.  
Of course, in the situation of Corollary~\ref{cor-Hopf-weak},
Theorem~\ref{thm-qc-ess-surj} gives us an equivalence of categories
between $(A,\Gamma )$-comodules and $(B, \Gamma _{f})$-comodules.  This
equivalence of categories takes an $(A, \Gamma )$-comodule $M$ to
$B\otimes _{A}M$.  

\section{Formal groups}\label{sec-formal}

In this section, we apply Corollary~\ref{cor-Hopf-weak} and the theory
of formal group laws to prove Theorem~\ref{main-E}.  We also recover the
change of rings theorems of Miller-Ravenel~\cite{miller-ravenel} and
Hovey-Sadofsky~\cite{hovey-sadofsky-picard}.  

This section requires familiarity with formal group laws and how they
are used in algebraic topology.  A good source for this material
is~\cite{ravenel}, especially Appendix~2 for formal group laws and
Chapter~4 for their use in algebraic topology.  

Fix a prime $p$ for use throughout this section.  Recall that
$(BP_{*},BP_{*}BP)$ is the universal Hopf algebroid for $p$-typical
formal group laws.  Here $BP_{*}=\Z _{(p)}[v_{1},v_{2},\dots ]$, and
$BP_{*}BP=BP_{*}[t_{1},t_{2},\dots ]$; see~\cite[Section~4.1]{ravenel}.
The fact that $(BP_{*},BP_{*}BP)$ is universal means that a $p$-typical
formal group law over a ring $R$ is equivalent to a ring homomorphism
$BP_{*}\xrightarrow{}R$, and a strict isomorphism of $p$-typical formal
group laws over $R$ is equivalent to a ring homomorphism
$BP_{*}BP\xrightarrow{}R$.  In case $R$ is graded, let us call a
$p$-typical formal group law over $R$ \emph{homogeneous} if its
classifying map $BP_{*}\xrightarrow{}R$ preserves the grading.  (An
example of a non-homogeneous formal group law is the formal group law
over $\Fp $ whose classifying map takes $v_{i}$ to $0$ for $i\neq n$ and
$v_{n}$ to $1$).  

Recall also the invariant ideal $I_{n}=(p,v_{1},\dots ,v_{n-1})$.  The
element $v_{n}$ is a primitive modulo $I_{n}$.  This means that there is
a Hopf algebroid 
\[
(A, \Gamma )=(v_{n}^{-1}BP_{*}/I_{n}, v_{n}^{-1}BP_{*}BP/I_{n}).
\]

\begin{definition}\label{defn-strict-height}
A $p$-typical formal group law over a ring $R$ is said to have
\emph{strict height $n$} if its classifying map factors through
$v_{n}^{-1}BP_{*}/I_{n}$.  
\end{definition}

Our application of Theorem~\ref{thm-qc-ess-surj} is then the following
theorem, which is Theorem~\ref{main-E} of the introduction.

\begin{theorem}\label{thm-awesome}
Fix a prime $p$ and an integer $n>0$.  Let $(A,\Gamma )$ denote the Hopf
algebroid $(v_{n}^{-1}BP_{*}/I_{n}, v_{n}^{-1}BP_{*}BP/I_{n})$.  Suppose
$B$ is a graded ring equipped with a homogeneous $p$-typical formal
group law of strict height $n$, classified by $f\mathcolon
A\xrightarrow{}B$.  Then the functor that takes an $(A,\Gamma
)$-comodule $M$ to $B\otimes _{A}M$ defines an equivalence of categories
from graded $(A,\Gamma )$-comodules to graded $(B,\Gamma _{f})$-comodules.
\end{theorem}

\begin{proof}
Let $D=A\otimes _{\Fp [v_{n},v_{n}^{-1}]}B$.  Let $x\mathcolon
A\xrightarrow{}D$ denote the ring homomorphism defined by $x(a)=a\otimes
1$, and let $y\mathcolon B\xrightarrow{}D$ denote the ring homomorphism
defined by $y(b)=1\otimes b$.  Then $x$ and the composite $yf$ induce
two formal group laws $F$ and $G$ over $D$, both $p$-typical and of
strict height $n$.  Furthermore, $x(v_{n})=yf(v_{n})$.  A result of
Lazard, as modified by
Strickland~\cite[Theorem~3.4]{hovey-sadofsky-picard}, then implies that
there is a faithfully flat graded ring extension $h\mathcolon
D\xrightarrow{}C$ and a strict isomorphism from $h_{*}G$ to $h_{*}F$.
This strict isomorphism is represented by a graded ring homomorphism $\alpha
\mathcolon \Gamma \xrightarrow{}C$.  Let $g\mathcolon B\xrightarrow{}C$
be the composite $hy$.  Since the domain of $\alpha $ is $h_{*}G$,
$\alpha \eta _{L}=gf\mathcolon A\xrightarrow{}C$.  This means that there
is a well-defined map
\[
g\otimes \alpha \mathcolon B \otimes _{A}\Gamma \xrightarrow{g\otimes
\alpha } C_{gf} \otimes _{A} \bimod{C}{\alpha \eta _{L}}{\ }
\xrightarrow{\mu } C.  
\]
Furthermore, $(g\otimes \alpha )\circ (f\otimes \eta _{R})$ represents
the codomain of $\alpha $, so is $hx$.  We know already that $h$ is a
faithfully flat ring extension, and we claim that $x$ is also a
faithfully flat ring extension.  Indeed, since $\Fp [v_{n},v_{n}^{-1}]$
is a graded field, $B$ is a free $\Fp [v_{n},v_{n}^{-1}]$-module, and so
$x$ makes $D$ into a free $A$-module.  Corollary~\ref{cor-Hopf-weak} and
Theorem~\ref{thm-qc-ess-surj} complete the proof. 
\end{proof}

In particular, we can take $B=E(m)_{*}/I_{n}$, where $m\geq n$ and
$E(m)$ is the Landweber exact Johnson-Wilson homology theory introduced
in~\cite{johnson-wilson}.  This leads to the following corollary. 

\begin{corollary}\label{cor-awesome-2}
Let $p$ be a prime and $m\geq n>0$ be integers.  Then the functor that takes
$M$ to $E(m)_{*}\otimes _{BP_{*}}M$ defines an equivalence of categories
\begin{multline*}
(v_{n}^{-1}BP_{*}/I_{n},v_{n}^{-1}BP_{*}BP/I_{n})\text{-comodules} \\
\xrightarrow{}
(v_{n}^{-1}E(m)_{*}/I_{n},v_{n}^{-1}E(m)_{*}E(m)/I_{n})\text{-comodules}.  
\end{multline*}
\end{corollary}

Using the method of~\cite{miller-ravenel}, we then get the following 
change of rings theorem, which is Theorem~\ref{main-F} of the
introduction. 

\begin{theorem}\label{thm-change-of-rings}
Let $p$ be a prime and $m\geq n>0$ be integers.  Suppose $M$ and $N$ are
$BP_{*}BP$-comodules such that $v_{n}$ acts isomorphically on $N$.  If
either $M$ is finitely presented, or if $N=v_{n}^{-1}N'$ where $N'$ is
finitely presented and $I_{n}$-nilpotent, then
\[
\text{Ext}^{**}_{BP_{*}BP}(M,N) \cong
\text{Ext}^{**}_{E(m)_{*}E(m)}(E(m)_{*}\otimes
_{BP_{*}}M,E(m)_{*}\otimes _{BP_{*}}N).  
\]
\end{theorem}

Note that, when $M=BP_{*}$, this is the Hovey-Sadofsky change of rings
theorem~\cite[Theorem~3.1]{hovey-sadofsky-picard}.  When $m=n$ and
$M=BP_{*}$, we get the Miller-Ravenel change of rings
theorem~\cite[Theorem~3.10]{miller-ravenel}.  

\begin{proof}
By Lemma~3.11 of~\cite{miller-ravenel}, $N$ is the direct limit of
comodules $v_{n}^{-1}N'$, where $N'$ is finitely presented and
$I_{n}$-nilpotent.  Since we are assuming either that $M$ is finitely
presented or that $N=v_{n}^{-1}N'$, in either case we may as well take
$N=v_{n}^{-1}N'$.  Then Lemma~3.12 of~\cite{miller-ravenel} reduces us
to the case $N=v_{n}^{-1}BP_{*}/I_{n}$.  In this case, one can check
using the cobar resolution (as
in~\cite[Proposition~1.3]{miller-ravenel}) that we have canonical
isomorphisms 
\[
\text{Ext}^{**}_{BP_{*}BP} (M,N) \cong
\text{Ext}^{**}_{v_{n}^{-1}BP_{*}BP/I_{n}}(v_{n}^{-1}M/I_{n},N) 
\]
and
\begin{multline*}
\text{Ext}^{**}_{E(m)_{*}E(m)} (E(m)_{*}\otimes _{BP_{*}}
M,E(m)_{*}\otimes _{BP_{*}} N) \cong \\
\text{Ext}^{**}_{v_{n}^{-1}E(m)_{*}E(m)/I_{n}} (E(m)_{*}\otimes _{BP_{*}} 
v_{n}^{-1}M/I_{n},E(m)_{*}\otimes _{BP_{*}} N)  
\end{multline*}
Now Corollary~\ref{cor-awesome-2} implies that 
\begin{multline*}
\text{Ext}^{**}_{v_{n}^{-1}BP_{*}BP/I_{n}}(v_{n}^{-1}M/I_{n},N) \cong \\
\text{Ext}^{**}_{v_{n}^{-1}E(m)_{*}E(m)/I_{n}} (E(m)_{*}\otimes _{BP_{*}} 
v_{n}^{-1}M/I_{n},E(m)_{*}\otimes _{BP_{*}} N).
\end{multline*}
This completes the proof.  
\end{proof}


\providecommand{\bysame}{\leavevmode\hbox to3em{\hrulefill}\thinspace}

\end{document}